\documentclass[reqno, draft]{amsart}
\usepackage{amsmath,amsthm,amscd,amssymb,latexsym,upref}

\newcommand{\bbN}{{\mathbb{N}}}
\newcommand{\bbR}{{\mathbb{R}}}
\newcommand{\bbP}{{\mathbb{P}}}
\newcommand{\bbZ}{{\mathbb{Z}}}
\newcommand{\bbC}{{\mathbb{C}}}

\newcommand{\bbQ}{{\mathbb{Q}}}

\newcommand{\calB}{{\mathcal B}}
\newcommand{\cB}{{\mathcal B}}
\newcommand{\calC}{{\mathcal C}}
\newcommand{\calD}{{\mathcal D}}
\newcommand{\calE}{{\mathcal E}}
\newcommand{\calF}{{\mathcal F}}

\newcommand{\cH}{{\mathcal H}}
\newcommand{\calK}{{\mathcal K}}
\newcommand{\calL}{{\mathcal L}}
\newcommand{\calM}{{\mathcal M}}

\newcommand{\calP}{{\mathcal P}}

\newcommand{\calV}{{\mathcal V}}

\renewcommand{\gg}{p}

\newcommand{\hatt}{\widehat}  

\newcommand{\Div}{\operatorname{Div}}

\newcommand{\no}{\nonumber}
\newcommand{\lb}{\label}
\newcommand{\f}{\frac}
\newcommand{\ul}{\underline}
\newcommand{\ol}{\overline}
\newcommand{\ti}{\tilde}
\newcommand{\wti}{\widetilde}

\newcommand{\uc}{{\underline{c}}}

\newcommand{\Oh}{O}
\newcommand{\oh}{o}

\newcommand{\ran}{\text{\rm{ran}}}

\newcommand{\dom}{\text{\rm{dom}}}
\newcommand{\bi}{\bibitem}
\newcommand{\humu}{{ \hat{\underline{\mu} }}}

\newcommand{\hmu}{{\hat{\mu} }}

\newcommand{\uz}{{\underline{z}}}

\newcommand{\uxi}{{\underline{\Xi}}}
\newcommand{\ome}{\omega}

\newcommand{\al}{\alpha}
\newcommand{\ual}{{\underline{\alpha}}}
\newcommand{\ua}{{\underline{A}}}
\newcommand{\hua}{{ \underline{\hatt{A} }}}

\newcommand{\uU}{{\underline{U}}}

\newcommand{\Sp}{\text{\rm Sp}}

\renewcommand{\Re}{\text{\rm Re}}
\renewcommand{\Im}{\text{\rm Im}}

\DeclareMathOperator{\diag}{diag}

\DeclareMathOperator{\sTl}{s-Tl}

\DeclareMathOperator{\sym}{Sym}
\DeclareMathOperator{\AP}{AP}
\DeclareMathOperator{\QP}{QP}

\newcommand{\symgg}{{\sym^\gg (\calK_\gg)}}

\allowdisplaybreaks 
\numberwithin{equation}{section}

\newtheorem{theorem}{Theorem}[section]
\newtheorem{lemma}[theorem]{Lemma}

\theoremstyle{definition}
\newtheorem{definition}[theorem]{Definition}
\newtheorem{hypothesis}[theorem]{Hypothesis}
\newtheorem{remark}[theorem]{Remark}

\begin{document}
\title[On the spectrum of Jacobi operators]{On the spectrum of Jacobi
operators with quasi-periodic algebro-geometric coefficients}
\author[V.\ Batchenko]{Vladimir Batchenko}
\address{Department of Mathematics,
University of Missouri, Columbia, MO 65211, USA}
\email{batchenv@math.missouri.edu}
\author[F.\ Gesztesy]{Fritz Gesztesy}
\address{Department of Mathematics,
University of Missouri, Columbia, MO 65211, USA}
\email{fritz@math.missouri.edu}
\urladdr{http://www.math.missouri.edu/personnel/faculty/gesztesyf.html}
\date{June 8}
\subjclass[2000]{Primary 34L05, 37K10, 37K20, 47B36; Secondary 34K14, 35Q51,
35Q58.} \keywords{Toda hierarchy, Jacobi operator, spectral theory.}
\thanks{To appear in {\it Intl. Math. Res. Papers (IMRP)}}

\begin{abstract}
We characterize the spectrum of one-dimensional Jacobi operators
$H=aS^{+}+a^{-}S^{-}+b$ in $l^{2}(\bbZ)$ with quasi-periodic
complex-valued algebro-geometric coefficients (which satisfy one
(and hence infinitely many) equation(s) of the stationary Toda
hierarchy) associated with nonsingular hyperelliptic curves. The
spectrum of $H$ coincides with the conditional stability set of
$H$ and can explicitly be described in terms of the mean value of
the Green's function of $H$.

As a result, the spectrum of $H$ consists of finitely many simple
analytic arcs in the complex plane. Crossings as well as
confluences of spectral arcs are possible and discussed as well.
\end{abstract}

\maketitle

\section{Introduction}\lb{s1}

It is well-known from the work of Date and Tanaka \cite{DT176},
\cite{DT276}, Dubrovin, Matveev, and Novikov \cite{DMN76},
Flaschka \cite{Fl175}, McKean \cite{McK79}, 
McKean and van Moerbeke \cite{MM75}, \cite{MM80}, Mumford \cite{Mu77},
Novikov, Manakov, Pitaevskii, and Zakharov \cite{NMPZ84}, Teschl \cite[Chs.
9,13]{Te00}, Toda \cite[Ch. 4]{To89}, \cite[Chs. 26-30]{To189}, 
van Moerbeke \cite{Mo79}, van Moerbeke and Mumford \cite{MM79},
that the self-adjoint Jacobi operator
\begin{equation}
H=aS^{+}+a^{-}S^{-}+b,\quad \text{dom}(H)=\ell^2(\bbZ), \lb{1.1}
\end{equation}
in $\ell^2(\bbZ)$ with real-valued periodic, or more generally,
algebro-geometric
{\it quasi-periodic} and {\it real-valued} coefficients $a$ and $b$ (i.e.,
coefficients that satisfy one (and hence infinitely many) equation(s) of the
stationary Toda hierarchy), leads to a finite-gap, or perhaps more
appropriately,
to a finite-band spectrum $\sigma (H)$ of the form
\begin{equation}
\sigma(H)=\bigcup_{m=1}^{p+1} [E_{2m-2},E_{2m-1}], \quad
E_0<E_1<\dots <E_{2p+1}. \lb{1.1aaa}
\end{equation}
Compared to the real-valued case, the corresponding spectral
properties of Jacobi operators with {\it periodic}
and {\it complex-valued} coefficients $a$ and $b$, to the best of our
knowledge, have been studied rather sparingly in the literature. We are
only aware of two papers by Na{\u\i}man \cite{Na62}, \cite{Na64}, in which
it is shown that the spectrum consists of a set of piecewise analytic arcs
which may have common endpoints. (This is in anology to the case of
(non-self-adjoint) one-dimensional periodic Schr\"odinger operators, cf.\
\cite{Ro63}, \cite{Se60}.) It seems plausible that the latter case is
connected with (complex-valued) stationary solutions of equations of the
Toda hierarchy. In particular, the next scenario in line, the determination
of the spectrum of $H$ in the case of {\it quasi-periodic} and {\it
complex-valued} solutions of the stationary Toda hierarchy apparently has
never been clarified. The latter spectral problems have been  open since
the late seventies and it is the purpose of this paper to provide a
comprehensive solution of them. For theta function representations of $a$
and $b$ in the complex-valued algebro-geometric case (without addressing
spectral theoretic questions) we refer, for instance, to Aptekarev
\cite{Ap86}, Dubrovin, Krichever, and Novikov
\cite{DKN90}, Krichever \cite{Kr78}--\cite{Kr83} (cf. also the appendix
written by Krichever in \cite{Du81}).

To describe our results, a bit of preparation is needed. Let
\begin{equation}
G(z,n,n')=(H-z)^{-1}(n,n'), \quad z\in\bbC\backslash\sigma(H), \;
n,n'\in\bbZ,  \lb{1.1a}
\end{equation}
be the Green's function of $H$ (here $\sigma(H)$ denotes the
spectrum of $H$) and denote by $g(z,n)$ the corresponding diagonal
Green's function of $H$ defined by
\begin{align}
&g(z,n)=G(z,n,n) =\f{\prod_{j=1}^\gg
[z-\mu_j(n)]}{R_{2\gg+2}(z)^{1/2}},
\lb{1.2a} \\
& R_{2\gg+2}(z)=\prod_{m=0}^{2\gg+1} (z-E_m), \quad
\{E_m\}_{m=0}^{2\gg+1}\subset\bbC,  \lb{1.2b} \\
& E_m\neq E_{m'}\, \text{  for $m\neq m'$, \;
$m,m'=0,1,\dots,2\gg+1$.} \lb{1.2c}
\end{align}
For any quasi-periodic (in fact, Bohr (uniformly) almost periodic)
sequence $f=\{f(k)\}_{k \in \bbZ}$ the mean value $\langle
f\rangle$ of $f$ is defined by
\begin{equation}
\langle f\rangle =\lim_{N\to\infty}\f{1}{2N+1} \sum_{k=-N}^{N}\,
f(k). \lb{1.3}
\end{equation}
Moreover, we introduce the set $\Sigma$ by
\begin{equation}
\Sigma=\bigg\{\lambda\in\bbC\,\bigg|\, \Re\bigg(\bigg\langle
\text{ln}\bigg(\frac{G_{\gg+1}(\lambda,\cdot)-y}{G_{\gg+1}(\lambda,\cdot)+y}\bigg)\bigg\rangle\bigg)=0\bigg\},
\lb{1.4}
\end{equation}
where $y^2=R_{2\gg+2}(z)$, and the polynomial $G_{\gg+1}(z,n)$ of degree
$p+1$ in $z$ will be defined in
\eqref{1.2.11b}. Here we observe that $G_{p+1}$ is given in terms of the
off-diagonal Green's function $G(z,n+1,n)$ by
\begin{equation}
\f{G_{p+1}(z,n)}{R_{2p+2}(z)^{1/2}}=1-2a(n)G(z,n+1,n)
\end{equation}
(cf.\ also \eqref{3.11}). In addition, we note that
\begin{equation}
\langle g(z,\cdot)\rangle=\f{\prod_{j=1}^\gg
\big(z-\wti\lambda_j\big)}{R_{2\gg+2}(z)^{1/2}} \lb{1.4a}
\end{equation}
for some constants $\{\wti\lambda_j\}_{j=1}^\gg\subset\bbC$.

Finally, we denote by $\sigma_{\rm p}(T)$, $\sigma_{\rm r}(T)$,
$\sigma_{\rm c}(T)$, $\sigma_{\rm{e}}(T)$, and
$\sigma_{\rm{ap}}(T)$, the point spectrum (i.e., the set of
eigenvalues), the residual spectrum, the continuous spectrum, the
essential spectrum (cf.\ \eqref{5.22b}), and the approximate point
spectrum of a densely defined closed operator $T$ in a complex
Hilbert space, respectively.

Our principal new results, to be proved in Section \ref{s4}, then
read as follows:

\begin{theorem}  \lb{t1.1}
Assume that $a$ and $b$ are quasi-periodic $($complex-valued\,$)$
solutions of the $\gg$th stationary Toda equation associated with
the hyperelliptic curve $y^2=R_{2\gg+2}(z)$ subject to
\eqref{1.2b} and \eqref{1.2c}. Then the following assertions hold:
\\
$(i)$ The point spectrum and residual spectrum of $H$ are empty
and hence the spectrum of $H$ is purely continuous,
\begin{align}
&\sigma_{\rm p}(H)=\sigma_{\rm r}(H)=\emptyset,  \lb{1.5} \\
&\sigma(H)=\sigma_{\rm c}(H)=\sigma_{\rm e}(H)=\sigma_{\rm ap}(H).
\lb{1.6}
\end{align}
$(ii)$ The spectrum of $H$ coincides with $\Sigma$ and equals the
conditional stability set of $H$,
\begin{align}
\sigma(H) &=\bigg\{\lambda\in\bbC\,\bigg|\, \Re\bigg(\bigg\langle
\ln\bigg(\frac{G_{\gg+1}(\lambda,\cdot)-y}
{G_{\gg+1}(\lambda,\cdot)+y}\bigg)\bigg\rangle\bigg)=0\bigg\}
\lb{1.10} \\
&=\{\lambda\in\bbC\,|\, \text{there exists at least one bounded
solution}  \no \\
& \hspace*{3.75cm} \text{$0\neq\psi\in \ell^\infty(\bbZ)$ of
$H\psi=\lambda\psi$}\}.  \lb{1.11}
\end{align}
$(iii)$ $\sigma(H)\subset \bbC$ is bounded,
\begin{equation}
\sigma(H)\subset \{z\in\bbC\,|\, \Re(z)\in [M_1,M_2],\, \Im(z)\in
[M_3,M_4]\}, \lb{1.12}
\end{equation}
where
\begin{align}
\begin{split}
& M_1=-2\sup_{n\in\bbZ}[|\Re (a(n))|]+\inf_{n\in\bbZ}[\Re(b(n))], \\
& M_2=2\sup_{n\in\bbZ}[|\Re (a(n))|]+\sup_{n\in\bbZ}[\Re(b(n))],  \lb{1.13}\\
& M_3=-2\sup_{n\in\bbZ}[|\Im (a(n))|]+\inf_{n\in\bbZ}[\Im(b(n))], \\
& M_4=2\sup_{n\in\bbZ}[|\Im(a(n))|]+\sup_{n\in\bbZ}[\Im(b(n))].
\end{split}
\end{align}
$(iv)$ $\sigma(H)$ consists of finitely many simple analytic arcs.
These analytic arcs may only end at the points
$\wti\lambda_1,\dots,\wti\lambda_\gg$, $E_0,\dots,E_{2\gg+1}$. \\
$(v)$ Each $E_m$, $m=0,\dots,2\gg+1$, is met by at least one of
these arcs. More precisely, a particular $E_{m_0}$ is hit by
precisely $2N_0+1$ analytic arcs, where $N_0\in\{0,\dots,\gg\}$
denotes the number of $\wti\lambda_j$ that coincide with
$E_{m_0}$. Adjacent arcs meet at an angle $2\pi/(2N_0+1)$ at
$E_{m_0}$. $($Thus,
generically, $N_0=0$ and precisely one arc hits $E_{m_0}$.$)$ \\
$(vi)$ Crossings of spectral arcs are permitted and take place
precisely when
\begin{align}
\begin{split}
&\Re\bigg(\bigg\langle
\ln\bigg(\frac{G_{\gg+1}(\wti\lambda_{j_0},\cdot)-y}{G_{\gg+1}(\wti\lambda_{j_0},\cdot)+y}\bigg)\bigg\rangle\bigg)=0
\lb{1.14} \\
&\quad \text{for some $j_0\in\{1,\dots,\gg\}$ with
$\wti\lambda_{j_0}\notin \{E_m\}_{m=0}^{2\gg+1}$}.
\end{split}
\end{align}
In this case $2M_0+2$ analytic arcs are converging toward
$\wti\lambda_{j_0}$, where $M_0\in\{1,\dots,\gg\}$ denotes the
number of $\wti\lambda_j$ that coincide with $\wti\lambda_{j_0}$.
Adjacent arcs meet at an angle $\pi/(M_0+1)$ at
$\wti\lambda_{j_0}$. $($Thus, if crossings occur, generically,
$M_0=1$ and two arcs cross at a right angle.$)$\\
$(vii)$ The resolvent set $\bbC\backslash\sigma (H)$ of $H$ is
path-connected.
\end{theorem}

Naturally, Theorem \ref{t1.1} applies to the special case where
$a$ and $b$ are periodic complex-valued solutions of the $\gg$th
stationary Toda equation associated with a nonsingular
hyperelliptic curve. Even in this special case, Theorem \ref{t1.1}
yields new facts which go beyond the previous results by Na{\u\i}man
\cite{Na62}, \cite{Na64}.

For analogous results in the context of one-dimensional Schr\"odinger
operators with quasi-periodic algebro-geometric KdV potentials we refer to
\cite{BG05}. 

Theorem \ref{t1.1} focuses on stationary quasi-periodic solutions of the
Toda hierarchy for the following reasons. First of all, the class of
algebro-geometric solutions of the (time-dependent) Toda hierarchy is
defined as the class of all solutions of some (and hence infinitely many)
equations of the stationary Toda hierarchy. Secondly, time-dependent
algebro-geometric solutions of a particular equation of the
(time-dependent) Toda hierarchy just represent isospectral deformations
(the deformation parameter being the time variable) of fixed stationary
algebro-geometric Toda solutions (the latter can be viewed as the initial
condition at a fixed time $t_0$). In the present case of quasi-periodic
algebro-geometric solutions of the $\gg$th Toda equation, the
isospectral manifold of such given solutions is a complex
$\gg$-dimensional torus, and time-dependent solutions trace out a
path in that isospectral torus (cf.\ the discussions in
\cite{GH03} and \cite{GH05}).

Finally, we give a brief discussion of the contents of each
section. In Section \ref{s2} we provide the necessary background
material including a quick construction of the Toda hierarchy of
nonlinear evolution equations and its Lax pairs using a polynomial
recursion formalism. We also discuss the hyperelliptic Riemann
surface underlying the stationary Toda hierarchy, the
corresponding Baker--Akhiezer function, and the necessary
ingredients to describe the analog of the Its--Matveev formula for
stationary Toda solutions. Section \ref{s3} focuses on the Green's
function of the Jacobi operator $H$, a key ingredient in our
characterization of the spectrum $\sigma(H)$ of $H$ in Section
\ref{s4} (cf.\ \eqref{1.10}). Our principal Section \ref{s4} is
then devoted to a proof of Theorem \ref{t1.1}. Appendix \ref{sA}
provides the necessary summary of tools needed from elementary
algebraic geometry (most notably the theory of compact
(hyperelliptic) Riemann surfaces) and sets the stage for some of
the notation used in Sections \ref{s2}--\ref{s4}. Appendix
\ref{sB} provides additional insight into one ingredient of the
theta function representation of the coefficients $a$ and $b$.

\section{The Toda Hierarchy, Hyperelliptic Curves, and Theta Function
Representations of the Coefficients $a$ and $b$} \label{s2}

In this section we briefly review the recursive construction of
the Toda hierarchy and associated Lax pairs following
\cite{BGHT98} and \cite{GH05}. Moreover, we discuss the
class of algebro-geometric solutions of the Toda hierarchy
corresponding to the underlying hyperelliptic curve and recall the
analog of the Its--Matveev formula for such solutions. The
material in this preparatory section is known and detailed
accounts with proofs can be found, for instance, in \cite{BGHT98}.
For the notation employed in connection with elementary concepts
in algebraic geometry (more precisely, the theory of compact
Riemann surfaces), we refer to Appendix \ref{sA}.

Throughout this section we assume that
\begin{equation}
a, b \in \ell^\infty (\bbZ), \quad a(n)\neq 0 \, \text{ for all 
$n\in\bbZ$},  \lb{1.2.1}
\end{equation}
and consider the second-order Jacobi difference expression
\begin{equation}
L=aS^{+}+a^{-}S^{-}+b, \lb{1.2.2}
\end{equation}
where $S^{\pm}$ denote the shift operators
\begin{equation}
(S^{\pm}f)(n)=f^{\pm}(n)=f(n{\pm}1), \quad n\in\bbZ, \; f \in
\ell^{\infty}(\bbZ). \lb{1.2.3}
\end{equation}

To construct the stationary Toda hierarchy we need a second
difference expression of order $2\gg +2, \, \gg \in {\bbN_{0}},$
defined recursively in the following. We take the quickest route
to the construction of $P_{2\gg +2}$, and hence to the Toda
hierarchy, by starting from the recursion relations
\eqref{1.2.4a}--\eqref{1.2.4c} below.

Define $\{{f_j}\}_{j \in \bbN_0}$ and $\{{g_j}\}_{j \in \bbN_0}$
recursively by
\begin{align}
& f_0=1, \quad g_0=-c_1, \lb{1.2.4a} \\
& 2f_{j+1}+g_j+g_j^{-}-2bf_j=0, \quad j\in \bbN_0, \lb{1.2.4b} \\
&
g_{j+1}-g_{j+1}^{-}+2\big(a^2f_{j}^{+}-(a^{-})^2f_j^{-}\big)-b(g_j-g_j^{-})=0,
\quad j \in \bbN_{0}. \lb{1.2.4c}
\end{align}
Explicitly, one finds
\begin{align}
& f_0=1, \no \\
& f_1=b+c_1, \no \\
& f_2=a^2+(a^{-})^2+b^2+c_1b+c_2, \, \text{ etc.,} \lb{1.2.5}\\
& g_0=-c_1, \no\\
& g_1=-2a^2-c_2, \no \\
& g_2=-2a^2(b^+ +b)+ c_1(-2 a^2)-c_3, \, \text{ etc.} \no
\end{align}
Here $\{ c_j\}_{j \in \bbN}$ denote undetermined summation
constants which naturally arise when solving
\eqref{1.2.4a}--\eqref{1.2.4c}.

Subsequently, it will be convenient to also introduce the
corresponding homogeneous coefficients $\hat f_j$ and $\hat g_j$,
defined by vanishing of the constants $c_k,\, k \in \bbN$,
\begin{align}
& \hat f_0=1, \quad \hat f_j=f_j \big|_{c_k=0, \, k=1,\dots,j},
\quad j\in\bbN, \no \\
& \hat g_j=g_j \big|_{c_k=0, \, k=1,\dots,j+1}, \quad j\in\bbN_0.
\lb{1.2.6}
\end{align}
Hence,
\begin{equation}
f_j=\sum_{k=0}^{j}c_{j-k}\hat{f}_{k},\quad
g_j=\sum_{k=1}^{j}c_{j-k}\hat{g}_{k}-c_{j+1},\quad j\in\bbN_0,
\lb{1.2.7}
\end{equation}
introducing
\begin{equation}
c_0=1. \lb{1.2.7a}
\end{equation}

Next we define difference expressions $P_{2\gg+2}$ of order
$2\gg+2$ by
\begin{equation}
P_{2\gg+2}=-L^{\gg+1}+\sum_{j=0}^{\gg} \Big(
g_j+2af_jS^{+}\Big)L^{\gg-j}+f_{\gg+1},\quad \gg\in\bbN_0.
\lb{1.2.8}
\end{equation}
Using the recursion relations \eqref{1.2.4a}--\eqref{1.2.4c}, the
commutator of $P_{2\gg+2}$ and $L$ can be explicitly computed and
one obtains
\begin{align}
[P_{2\gg+2},L]=&-a\big(g_\gg^{+}+g_\gg+f_{\gg+1}^{+}+f_{\gg+1}-2b^{+}f_\gg^{+}
\big)S^{+} \no \\
&+2\big(-b(g_{\gg}+f_{\gg+1})+a^2f_\gg^{+}
-(a^{-})^2f_{\gg}^{-}+b^2f_\gg\big) \no \\
&-a^{-}\big(g_\gg+g_\gg^{-}+f_{\gg+1}+f_{\gg+1}^{-}-2bf_\gg
\big)S^{-}, \quad \gg\in\bbN_0. \lb{1.2.8a}
\end{align}
In particular, $(L,P_{2\gg+2})$ represents the celebrated
\textit{Lax pair} of the Toda hierarchy. Varying $\gg\in\bbN_0$,
the stationary Toda hierarchy is then defined in terms of the
vanishing of the commutator of $P_{2\gg+2}$ and $L$ in
\eqref{1.2.8a} by,
\begin{equation}
[P_{2\gg+2}, L] =\sTl_\gg(a,b)=0,\quad \gg\in\bbN_0. \lb{1.2.8aa}
\end{equation}
Thus one finds
\begin{align}
g_\gg+g_\gg^{-}+f_{\gg+1}+f_{\gg+1}^{-}-2bf_\gg&=0, \lb{1.2.8ab} \\
-b(g_\gg+f_{\gg+1})+a^2f_\gg^{+}-(a^{-})^2f_g^{-}+b^2f_\gg&=0.
\lb{1.2.8ac}
\end{align}
Using \eqref{1.2.4b} with $j=\gg$ one concludes that
\eqref{1.2.8ab} reduces to
\begin{equation}
f_{\gg+1}=f_{\gg+1}^{-}, \lb{1.2.8b}
\end{equation}
that is, $f_{\gg+1}$ is a lattice constant. Similarly, one infers
by subtracting $b$ times \eqref{1.2.8ab} from twice
\eqref{1.2.8ac} and using \eqref{1.2.4c} with $j=\gg$, that
$g_{\gg+1}$ is a lattice constant as well, that is,
\begin{equation}
g_{\gg+1}=g_{\gg+1}^{-}. \lb{1.2.8c}
\end{equation}
Equations \eqref{1.2.8b} and \eqref{1.2.8c} constitute the $\gg$th
stationary equation in the Toda hierarchy, which will be denoted
by
\begin{equation}
\sTl_\gg(a,b)=\begin{pmatrix} f^+_{p+1}-f_{p+1} \\
g_{p+1} -g_{p+1}^- \end{pmatrix}=0, \quad \gg\in\bbN_0.
\lb{1.2.9}
\end{equation}
Explicitly,
\begin{align}
& \sTl_0(a,b)=  \begin{pmatrix} b^+ -b \\
2\big((a^-)^2-a^2\big) \end{pmatrix} =0, \no \\
& \sTl_1(a,b)=  \begin{pmatrix}
     (a^+)^2-(a^-)^2+(b^+)^2-b^2  \\
2(a^-)^2(b+b^-) - 2a^2(b^{+}+b) \end{pmatrix} \\
& \quad\quad\quad\quad\quad\quad + c_1 
\begin{pmatrix}
b^+ -b \\
2 \big((a^-)^2 -a^2\big)
\end{pmatrix} =0, \quad \text{etc.,} \no
\end{align}
represent the first few equations of the stationary Toda
hierarchy. By definition, the set of solutions of \eqref{1.2.8aa},
with $\gg$ ranging in $\bbN_0$ and $c_k\in\bbC$, $k\in\bbN$,
represents the class of algebro-geometric Toda solutions.

In the following we will frequently assume that $a$ and $b$
satisfy the $\gg$th stationary Toda equation. By this we mean they
satisfy one of the $\gg$th stationary Toda equations after a
particular choice of integration constants $c_k\in\bbC$,
$k=1,\dots,\gg$, $\gg\in\bbN$, has been made.

Next, we introduce polynomials $F_\gg(z,n)$ and $G_{\gg+1}(z,n)$
of degree $\gg$ and $\gg+1$ with respect to the spectral parameter
$z\in\bbC$ by
\begin{align}
F_\gg(z,n)&=\sum_{j=0}^\gg z^jf_{\gg-j}(n), \lb{1.2.11a} \\
G_{\gg+1}(z,n)&=-z^{\gg+1}+\sum_{j=0}^\gg
z^jg_{\gg-j}(n)+f_{\gg+1}(n). \lb{1.2.11b}
\end{align}
Explicitly, one obtains
\begin{align}
F_0&=1, \no \\
F_1&=z+b+c_1, \no \\
F_2&=z^2+bz+a^2+(a^-)^2+b^2+c_1(z+b)+c_2, \, \text{ etc.,}
\lb{1.2.12}
\\ G_1&=-z+b,\no  \\
G_2&=-z^2+(a^-)^2-a^2+b^2+c_1(-z+b), \, \text{ etc.} \no 
\end{align}

Next, we study the restriction of the difference expression
$P_{2\gg+2}$ to the two-dimensional  kernel (i.e., the formal null
space in an algebraic sense as opposed to the functional analytic
one) of $(L-z)$. More precisely, let
\begin{align}
& \text{ker}(L-z)=\{\psi \colon \bbZ \to \bbC\cup\{\infty\} \mid
(L-z)\psi=0\}. \lb{1.2.13a}
\end{align}
Then \eqref{1.2.8} implies
\begin{equation}
P_{2\gg+2}\mid_{\text{ker}(L-z)}=\big(2aF_\gg(z)S^{+}+G_{\gg+1}(z)\big)\big|_{\text{ker}(L-z)}.
\lb{1.2.13}
\end{equation}
Therefore, the Lax relation \eqref{1.2.8a} becomes
\begin{align}
0&=[P_{2\gg+2},L]\mid_{\text{ker}(L-z)}=\Big(a\big(2(z-b^{+})F_\gg^{+}-2(z-b)F_\gg+G_{\gg+1}^{-}-G_{\gg+1}^{+}\big)S^{+}
\no \\
&\quad+\big(2(a^{-})^2F_\gg^{-}-2a^2F_\gg^{+}+(z-b)(G_{\gg+1}^{-}-G_{\gg+1})\big)
\Big)\Big|_{\text{ker}(L-z)}, \lb{1.2.14}
\end{align}
or, equivalently,
\begin{align}
2(z-b^{+})F_\gg^{+}-2(z-b)F_\gg+G_{\gg+1}^{+}-G_{\gg+1}^{-}&=0,
\lb{1.2.15a} \\
2a^2F_\gg^{+} - 2(a^{-})^2F_\gg^{-} +(z-b)(G_{\gg+1}-G_{\gg+1}^-)&=0.
\lb{1.2.15b}
\end{align}
Upon summing \eqref{1.2.15a} one infers
\begin{equation}
2(z-b^{+})F_\gg^{+}+G_{\gg+1}^{+}+G_{\gg+1}=0, \quad \gg\in
\bbN_0, \lb{1.2.16a}
\end{equation}
and inserting \eqref{1.2.15a} into \eqref{1.2.15b} then implies
\begin{equation}
(z-b)^2F_\gg+(z-b)G_{\gg+1}+a^2F_\gg^{+}-(a^{-})^2F_\gg^{-}=0,
\quad \gg\in \bbN_0. \lb{1.2.16b}
\end{equation}
Combining equations \eqref{1.2.15b} and \eqref{1.2.16a} one
concludes that the quantity
\begin{equation}
R_{2\gg+2}(z) = G_{\gg+1}(z,n)^2-4a(n)^2F_\gg(z,n)F_\gg^{+}(z,n)
\lb{1.2.17}
\end{equation}
is a lattice constant, and hence depends on $z$ only. Thus, one can
write
\begin{equation}
R_{2\gg+2}(z) = \prod_{m=0}^{2\gg+1} (z-E_m), \quad \{
E_m\}_{m=0}^{2\gg+1} \subset \bbC. \lb{1.2.18}
\end{equation}
One can show that equation \eqref{1.2.17} leads to an explicit
determination of the integration constants $c_1,\dots,c_\gg$ in
\begin{equation}
\sTl_\gg(a,b)=0 \lb{1.2.12A}
\end{equation}
in terms of the zeros $E_0,\dots,E_{2\gg+1}$ of the associated
polynomial $R_{2\gg+2}$ in \eqref{1.2.18}. In fact, one can prove that 
\begin{equation}
c_k=c_k(\ul E), \quad k=1,\dots,\gg, \lb{1.2.12B}
\end{equation}
where
\begin{align}
c_k(\ul E) \no &=-\!\!\!\!\!\!\!\sum_{\substack{j_0,\dots,j_{2\gg+1}=0\\
        j_0+\cdots+j_{2\gg+1}=k}}^{k}\!\!
\f{(2j_0)!\cdots(2j_{2\gg+1})!} {2^{2k} (j_0!)^2\cdots
(j_{2\gg+1}!)^2 (2j_0-1)\cdots(2j_{2\gg+1}-1)}
E_0^{j_0}\cdots E_{2\gg+1}^{j_{2\gg+1}}, \no \\
& \hspace*{8cm} k=1,\dots,p \label{1.2.12C}
\end{align}
are symmetric functions of $\ul E=(E_0,\dots,E_{2p+1})$.

\begin{remark} \lb{r2.1}
Since by \eqref{1.2.4b}, \eqref{1.2.4c}, \eqref{1.2.11a} and
\eqref{1.2.11b}, $a$ enters quadratically in $F_\gg$ and
$G_{\gg+1}$, the Toda hierarchy \eqref{1.2.9} is invariant under
the substitution
\begin{equation}
a(n)\rightarrow a_\epsilon(n)=\epsilon(n)a(n), \quad
\epsilon(n)\in\{+1,-1\}, \quad n \in \bbZ. \lb{1.2.12a}
\end{equation}
\end{remark}

We emphasize that the result \eqref{1.2.13} is valid independently
of whether or not $P_{2\gg+2}$ and $L$ commute.  However, the fact
that the two difference expressions $P_{2\gg+2}$ and $L$ commute
implies the existence of an algebraic relationship between them.
This gives rise to the Burchnall--Chaundy polynomial for the Toda
hierarchy first discussed in the discrete context by Na{\u\i}man
\cite{Na62}, \cite{Na64}.

\begin{theorem} \label{t1.2.1}
Fix $\gg\in\bbN_0$ and assume that $P_{2\gg+2}$ and $L$ commute,
$[P_{2\gg+2},L]=0$, or equivalently, suppose $\sTl_{\gg}(a,b)=0$.
Then $L$ and $P_{2\gg+2}$ satisfy an algebraic relationship of the
type $($cf.\ \eqref{1.2.18}$)$
\begin{align}
\begin{split}
&\calF_\gg(L,P_{2\gg+2})= P_{2\gg+2}^2-R_{2\gg+2}(L)=0, \label{1.2.19} \\
& R_{2\gg+2}(z) = \prod_{m=0}^{2\gg+1} (z-E_m), \quad z\in \bbC. 
\end{split}
\end{align}
\end{theorem}

The expression $\calF_\gg(L,P_{2\gg+2})$ is called the
Burchnall--Chaundy polynomial of the Lax pair $(L,P_{2\gg+2})$.
Equation \eqref{1.2.19} naturally leads to the hyperelliptic
curve $\calK_\gg$ of (arithmetic) genus $\gg\in\bbN_0$, where
\begin{align}
\begin{split}
&\calK_\gg \colon \calF_\gg(z,y)=y^2-R_{2\gg+2}(z)=0,  \lb{1.2.20} \\
&R_{2\gg+2}(z) = \prod_{m=0}^{2\gg+1} (z-E_m), \quad \{
E_m\}_{m=0}^{2\gg+1} \subset \bbC. 
\end{split}
\end{align}
The curve $\calK_\gg$ is compactified by joining two points
$P_{\infty_{\pm}}$, $P_{\infty_+}\neq P_{\infty_-},$ at
infinity. For notational simplicity, the resulting curve is still
denoted by $\calK_\gg$. Points $P$ on $\calK_\gg \backslash
{P_{\infty_{\pm}}}$ are  represented as pairs $P=(z,y)$, where
$y(\cdot)$ is the meromorphic function on $\calK_\gg$ satisfying
$\calF_\gg(z,y)=0$. The complex structure on $\calK_\gg$ is then
defined in the usual way, see Appendix \ref{sA} for the case of
nonsingular curves. Hence,
$\calK_\gg$ becomes a two-sheeted hyperelliptic Riemann surface of
(arithmetic) genus $\gg\in\bbN_0$ (possibly with a singular affine
part) in a standard manner.

We also emphasize that by fixing the curve $\calK_\gg$ (i.e., by
fixing $E_0,\dots,E_{2\gg+1}$), the integration constants
$c_1,\dots,c_\gg$ in the corresponding stationary $\sTl_\gg$
equation are uniquely determined as is clear from \eqref{1.2.12B}
and \eqref{1.2.12C}, which establish the integration constants
$c_k$ as symmetric functions of $E_0,\dots,E_{2\gg+1}$.

For notational simplicity we will usually tacitly assume that
$\gg\in\bbN$. The trivial case $\gg=0$, which leads to
$a(n)^2=(E_1-E_0)^2/16$, $b(n)=(E_0+E_1)/2$, $n\in\bbN$, is of no interest
to us in this paper.

In the following, the zeros\footnote{If $a,b\in
\ell^\infty(\bbZ)$, these zeros are the Dirichlet  eigenvalues of
a bounded operator on $\ell^2(\bbZ)$ associated with the difference
expression $L$ and a Dirichlet boundary condition  at $n\in\bbZ$.}
of the polynomial $F_\gg(\cdot,n)$ (cf.\ \eqref{1.2.11a}) will
play a special role. We denote them by $\{\mu_j(n)\}_{j=1}^\gg$
and write
\begin{equation}
F_\gg(z,n) =\prod_{j=1}^\gg [z-\mu_j(n)]. \lb{1.2.21}
\end{equation}
The next step is crucial; it permits us to ``lift'' the zeros
$\mu_j$ of $F_\gg$ from $\bbC$ to the curve $\calK_\gg$. {}From
\eqref{1.2.17} and \eqref{1.2.21} one infers
\begin{equation}
R_{2\gg+2}(z) - G_{\gg+1}(z)^2 = 0, \quad
z\in\{\mu_j, \mu_k^+\}_{j,k=1,\dots,\gg}. \lb{1.3.7}
\end{equation}
We now introduce $\{ \hat \mu_j(n) \}_{j=1,\dots,\gg}\subset
\calK_\gg$ by
\begin{align}
& \hat{\mu}_j(n)=\big( \mu_j(n),-G_{\gg+1}(\mu_j(n),n) \big),
\quad j=1,...,\gg, \; n\in\bbZ. \lb{1.2.24a}
\end{align}

Next, we recall equation \eqref{1.2.17} and define the
fundamental meromorphic function $\phi(\cdot,n)$ on $\calK_\gg$ by
\begin{align}
\phi(P,n)&=\frac{y-G_{\gg+1}(z,n)}{2a(n)F_\gg (z, n)} \lb{1.2.22a}\\
& = \frac{-2a(n)F_{\gg} (z,n+1)}{y+G_{\gg+1}(z,n)}, \lb{1.2.22b} \\
& P = (z,y)\in\calK_\gg, \; n\in\bbZ \no
\end{align}
with divisor $(\phi(\cdot, n))$ of $\phi(\cdot, n)$ given by
\begin{equation}
\big( \phi(\cdot,n) \big)=\calD_{P_{\infty_+} \hat{\underline{
\mu}}(n+1)}-\calD_{P_{\infty_-}\hat{\underline{\mu}}(n)},
\lb{1.2.23}
\end{equation}
using \eqref{1.2.21} and \eqref{1.2.24a}. Here we abbreviated
\begin{equation}
\hat{\underline \mu} = \{\hat \mu_1, \dots, \hat \mu_\gg\} \in
\symgg \lb{1.2.24b}
\end{equation}
(cf.\ the notation introduced in Appendix \ref{sA}). The
stationary Baker--Akhiezer function $\psi(\cdot,n,n_0)$ on
$\calK_\gg$ is then defined in terms of $\phi(\cdot,n)$ by
\begin{align}
& \psi(P,n,n_0)= \begin{cases} \prod_{m=n_0}^{n-1}\phi(P,m) &
\text{for}\quad n\geq n_0+1, \\
1 & \text{for}\quad n=n_0, \\
\prod_{m=n}^{n_0-1} {\phi(P,m)}^{-1} \quad &\text{for}\quad n\leq
n_0-1
\end{cases} \lb{1.2.25a}
\end{align}
with divisor $\big(\psi(\cdot,n,n_0)\big)$ of $\psi(P,n,n_0)$
given by
\begin{equation}
\big(\psi(\cdot,n,n_0)\big)=\calD_{\underline{\hat{\mu}}(n)}-\calD_{\underline{\hat{\mu}}(n_0)}+(n-n_0)(\calD_{P_{\infty
_ +}}-\calD_{P_{\infty_ -}}). \lb{1.2.25b}
\end{equation}
Basic properties of $\phi$ and $\psi$ are summarized in the
following result. We denote by  $W(f,g)(n)=a(f g^+ -f^+ g)$ the
Wronskian of two complex-valued sequences $f$ and $g$, and denote 
$P^*=(z,-y)$ for $P=(z,y)\in\calK_p$.

\begin{lemma} \lb{l1.3.1}
Assume \eqref{1.2.1} and suppose $a, b$ satisfy the
$\gg$th stationary Toda equation \eqref{1.2.9}. Moreover, let $P= (z,y)
\in \calK_\gg \backslash \{P_{\infty_{\pm}}\}$ and $(n,n_0) \in
\bbZ^2$.  Then $\phi$ satisfies the Riccati-type equation
\begin{align}
& a\phi(P) + a^{-}\phi^{-}(P)^{-1} = z-b, \lb{1.2.26}
\intertext{as well as}
& \phi(P) \phi (P^*)= \f{F^{+}_{\gg}(z)}{F_\gg (z)},\lb{1.2.27}\\
&\phi (P) + \phi (P^*)= \f{-G_{\gg+1} (z)}{aF_\gg (z)},\lb{1.2.28}\\
& \phi(P)-\phi (P^*)=\f{y(P)}{aF_\gg (z)}. \lb{1.2.29}
\end{align}
Moreover, $\psi$ satisfies
\begin{align}
& \big(L-z(P)\big) \psi(P) =0, \quad \big(P_{2\gg+2}-y(P)\big)
\psi(P) =0, \lb{1.2.30} \\
& \psi (P,n,n_0) \psi (P^*, n, n_0)= \f{F_\gg
(z,n)}{F_\gg(z,n_0)},\lb{1.2.31}\\
& a(n)\big[\psi(P,n,n_0)\psi(P^*,n+1,n_0)
+\psi(P^*,n,n_0)\psi(P,n+1,n_0)\big]\no\\
&\quad=-\f{G_{\gg+1}(z,n)}{F_\gg(z,n_0)},
\lb{1.2.32a} \\
&W(\psi (P,\cdot,n_0), \psi (P^*, \cdot,
n_0))=-\f{y(P)}{F_{\gg}(z,n_0)}. \lb{1.2.32ab}
\end{align}
\end{lemma}

Combining the polynomial recursion approach with \eqref{1.2.21}
readily yields trace formulas for the Toda invariants, which are
polynomial expressions of $a$ and $b$, in terms of the zeros $\mu_j$ of
$F_\gg$.

\begin{lemma} \lb{l1.3.9a}
Assume \eqref{1.2.1} and suppose $a, b$ satisfy the
$\gg$th stationary Toda equation \eqref{1.2.9}. Then,
\begin{align}
& a(n)^2=\f{1}{2}\sum_{j=1}^{\gg}
R_{2\gg+2}^{1/2}(\hat{\mu}_j(n))\prod_{\substack{k=1\\
k\neq j}}^\gg [\mu_j(n)-\mu_k(n)]^{-1}+\f{1}{4}[b^{(2)}(n) -
b(n)^2], \lb{1.2.21a} \\
& b(n)=\f{1}{2}\sum_{m=0}^{2\gg+1} E_m - \sum_{j=1}^\gg \mu_j(n),
\lb{1.2.21b} \\
& b^{(k)}(n)=\f{1}{2}\sum_{m=0}^{2\gg+1} E_m^k - \sum_{j=1}^\gg
\mu_j(n)^k,\quad k\in \bbN.
\end{align}
\end{lemma}

Strictly speaking, \eqref{1.2.21a} is only valid in the case where for
all $n\in\bbZ$, $\mu_j(n)\neq \mu_k(n)$ for $j\neq k$. The case where some
of the $\mu_j$ coincide requires a limiting argument that will be omitted as
\eqref{1.2.21a} will play no further role in this paper. The details of this
limiting procedure can be found in \cite{GHT05}.

From this point on we assume that the affine part of $\calK_\gg$ is
nonsingular, that is,
\begin{equation}
E_m\neq E_{m'} \text{  for $m\neq m'$, \;
$m,m'=0,1,\dots,2\gg+1$}.\lb{1.3.52A}
\end{equation}

Since nonspecial divisors play a fundamental role in this context
we also recall the following fact.

\begin{lemma} \lb{l1.3.9ba}
Assume \eqref{1.2.1} and suppose $a, b$ satisfy the
$\gg$th stationary Toda equation \eqref{1.2.9}. In addition, assume that 
the affine part of $\calK_\gg$ is nonsingular. Let $\calD_{\humu}$,
$\humu=\{\hmu_1,\dots,\hmu_\gg\}\in\symgg$, be the Dirichlet divisor of
degree $\gg$ associated with $a$, $b$ defined according to
\eqref{1.2.24a}, that is,
\begin{equation}
\hmu_j(n)=\big(\mu_j(n),-G_{\gg+1}(\mu_j(n),n)\big), \quad
j=1,\dots,\gg, \; n\in\bbZ. \lb{1.3.59AA}
\end{equation}
Then $\calD_{\humu(n)}$ is nonspecial for all $n\in\bbZ$.
Moreover, there exists a constant $C_{\mu}>0$ such that
\begin{equation}
|\mu_j(n)|\leq C_{\mu}, \quad j=1,\dots,\gg, \; n\in\bbZ.
\lb{1.3C}
\end{equation}
\end{lemma}

We remark that if $a,b \in \ell^\infty(\bbZ)$ satisfy the $\gg$th
stationary Toda equation \eqref{1.2.9}, then automatically $a(n)\neq 0$
for all $n\in\bbZ$ (cf.\ \cite{GHT05}).

We continue with the theta function representation for $\psi$,
$a$, and $b$. For general background information and the notation
employed we refer to Appendix \ref{sA}.

Let $\theta$ denote the Riemann theta function associated with
$\calK_\gg$ (whose affine part is assumed to be nonsingular) and a
fixed homology basis $\{a_j,b_j\}_{j=1}^\gg$. Next, choosing a
base point $Q_0 \in \calB(\calK_\gg)$ in the set of branch points
of $\calK_\gg$, we recall that the Abel maps $\ul{A}_{Q_0}$ and
$\ual_{Q_0}$ are defined by \eqref{a42} and \eqref{a43}, and the
Riemann vector $\underline{\Xi}_{Q_0}$ is given by \eqref{a55}.
Then Abel's theorem (cf. \eqref{a53}) and \eqref{1.2.25b} yields
\begin{align}
\begin{split}
       \ual_{Q_0} (\calD_{\humu(n)}) &= \ual_{Q_0}
(\calD_{\humu(n_0)}) - \ul{A}_{P_{\infty_-}}(P_{\infty_+})(n-n_0) \lb{1.2.33}\\
& =\ual_{Q_0}(\calD_{\humu(n_0)}) -
2\ul{A}_{Q_0}(P_{\infty_+})(n-n_0).
\end{split}
\end{align}

Next, let $\omega_{P_{\infty_+},P_{\infty_-}}^{(3)}$ denote the
normalized differential of the third kind defined by
\begin{align}
&\omega_{P_{\infty_+},P_{\infty_-}}^{(3)}  = \f1{ y}
\prod_{j=1}^\gg(z -\lambda_j) d z \underset{\zeta\to
0}{=}\pm\big(\zeta^{-1}+\Oh(1)\big)d\zeta \text{ as $P\to
P_{\infty_\pm}$},  \lb{1.2.34} \\
& \hspace*{8.85cm} \zeta=1 /z, \no
\end{align}
where the constants $\lambda_j\in\bbC$, $j=1,\dots, \gg$, are
determined by employing the normalization
\begin{equation}
\int_{a_j}\ome_{P_{\infty+},P_{\infty-}}^{(3)}=0, \quad j=1,\dots,
\gg. \lb{1.2.34aa}
\end{equation}
One then infers
\begin{equation}
\int_{Q_0}^P \ome_{P_{\infty+},P_{\infty-}}^{(3)}
\underset{\zeta\to 0}{=}\pm\ln (\zeta) + e^{(3)}_0(Q_0)+\Oh(\zeta)
\text{ as $P\to P_\infty$} \lb{1.3.71b}
\end{equation}
for some constant $e^{(3)}_0(Q_0)\in\bbC$. The vector of
$b$-periods of $\omega_{P_{\infty_+},P_{\infty_-}}^{(3)}/(2\pi i)$
is denoted by
\begin{equation}
\ul{U}_{0}^{(3)}=\big({U}_{0,1}^{(3)},\dots,{U}_{0,\gg}^{(3)}\big),
\quad {U}_{0,j}^{(3)}=\f{1}{2\pi i}\int_{b_j}
\omega_{P_{\infty_+},{P_{\infty_-}}}^{(3)},\quad j=1,\dots,\gg.
\lb{1.2.35}
\end{equation}
Since $Q_0$ is a branch point, $Q_0 \in \calB(\calK_\gg)$, one concludes 
by \eqref{a27a} that
\begin{equation}
{U}_{0}^{(3)}=\ul{A}_{P_{\infty -}}({P_{\infty
+}})=2\ul{A}_{Q_{0}}({P_{\infty +}}). \lb{1.2.36}
\end{equation}
In the following it will be convenient to introduce the abbreviation
\begin{align}
&\uz(P,\ul Q) =\ul\Xi_{Q_0}-\ul A_{Q_0}(P)+
\ul\alpha_{Q_0}(\calD_{\ul Q}), \lb{1.2.37}\\
& P\in\calK_\gg, \; \ul Q=\{Q_1,\dots, Q_\gg\}
\in \sym^\gg \no (\calK_\gg).
\end{align}
We note that $\ul{z}(\cdot,\ul{Q})$ is independent of the choice
of base point $Q_0$.

The zeros and the poles of $\psi$ as recorded in \eqref{1.2.25b}
suggest consideration of the following expression involving
$\theta$-functions
     (cf. \eqref{a31})
\begin{equation}
\frac{\theta\big(\underline{z}(P,\humu (n))\big)
}{\theta\big(\underline{z}(P,\humu (n_0))\big)} \text{exp}\bigg(
\int_{Q_0}^P \omega_{P_{\infty_+},P_{\infty_-}}^{(3)}\bigg).
\lb{1.2.38}
\end{equation} Here we agree to use the same path of integration
from $Q_0$ to $P$ on $\calK_\gg$ in the Abel map
$\underline{\hat{A}}_{Q_0}(P)$ in $\underline{z}(P,\humu(n))$ and
in the integral of $\omega_{P_{\infty_+},P_{\infty_-}}^{(3)}$ in
the exponent of \eqref{1.2.38}. With this convention the
expression \eqref{1.2.38} is well-defined on $\calK_\gg$ (cf.\ Remark
\ref{raa26a}, however) and one concludes
\begin{equation}
\psi(P,n,n_0)=C(n,n_0)\frac{\theta\big(\underline{z}(P,\humu(n))\big)}{\theta\big(\underline{z}(P,\humu(n_0))\big)}\text{exp}\bigg((n-n_0)
\int_{Q_0}^P \omega_{P_{\infty_+},P_{\infty_-}}^{(3)}\bigg).
\lb{1.2.39}
\end{equation}
To determine $C(n,n_0)$ one can use \eqref{1.2.31} for
$P=P_{\infty_+}$ and $P^{*}=P_{\infty_-}$. Hence,
\begin{equation}
C(n,n_0)^2=\frac{\theta\big(\underline{z}(P_{\infty_+},\humu(n_0))\big)\theta\big(\underline{z}(P_{\infty_+},\humu(n_0-1))\big)}
{\theta\big(\underline{z}(P_{\infty_+},\humu(n))\big)\theta\big(\underline{z}(P_{\infty_+},\humu(n-1))\big)}.
\lb{1.2.40}
\end{equation}
Thus, one obtains the following well-known result.

\begin{theorem} \lb{t1.3.10}
Assume \eqref{1.2.1} and suppose $a, b$ satisfy the $\gg$th stationary Toda
equation \eqref{1.2.9}. In addition,
assume the affine part of $\calK_\gg$ to be nonsingular and let
$P\in \calK_\gg \backslash \{ P_{\infty_\pm}\}$ and
$n,n_0\in\bbZ$. Then $\calD_{\ul{\hat\mu}(n)}$ is nonspecial for
$n\in\bbZ$. Moreover,\footnote{To avoid multi-valued expressions
in formulas such as \eqref{1.2.41}, etc., we agree to always
choose the same path of integration  connecting $Q_0$ and $P$ and
refer to Remark \ref{raa26a} for additional tacitly assumed
conventions.}
\begin{equation}
\psi(P,n,n_0)=
C(n,n_0)\frac{\theta(\underline{z}(P,\humu(n)))}{\theta(\underline{z}(P,\humu(n_0)))}\exp
\bigg((n-n_0) \int_{Q_0}^P
\omega_{P_{\infty_+},P_{\infty_-}}^{(3)}\bigg), \lb{1.2.41}
\end{equation}
where
\begin{equation}
C(n,n_0)=\bigg[
\frac{\theta\big(\underline{z}(P_{\infty_+},\humu(n_0))\big)\theta\big(\underline{z}(P_{\infty_+},\humu(n_0-1))\big)}
{\theta\big(\underline{z}(P_{\infty_+},\humu(n))\big)\theta\big(\underline{z}(P_{\infty_+},\humu(n-1))\big)}\bigg]^{1/2},
\lb{1.2.42}
\end{equation}
with the linearizing property of the Abel map,
\begin{equation}
\ual_{Q_0} (\calD_{\humu(n)}) = \big( \ual_{Q_0}
(\calD_{\humu(n_0)}) - 2\ua_{Q_0}(P_{\infty_+})(n-n_0)\big) \pmod
{L_\gg}. \lb{1.2.43}
\end{equation}
The coefficients $a$ and $b$ are given by
\begin{align}
&
a(n)=\tilde{a}
\bigg[\frac{\theta\big(\underline{z}(P_{\infty_+},\humu(n-1))\big)
\theta\big(\underline{z}(P_{\infty_+},\humu(n+1))\big)}
{\theta\big(\underline{z}(P_{\infty_+},\humu(n))\big)^2}\bigg]^{1/2}, \quad
     n\in\bbZ, \lb{1.2.44} \\
& b(n)=\frac{1}{2}\sum_{m=0}^{2\gg+1}E_m -\sum_{j=1}^\gg\lambda_j
+\sum_{j=1}^{\gg}c_j(\gg)\frac{\partial}{\partial\omega_j}\ln
\bigg[\frac{\theta\big(\underline{\omega}
+\underline{z}(P_{\infty_+},\humu(n))\big)}{\theta\big(\underline{\omega}+
\underline{z}(P_{\infty_+},\humu(n-1))\big)}\bigg]
\bigg|_{\underline{\omega}=0}, \no \\
& \hspace*{10cm} n\in\bbZ,  \lb{1.2.45}
\end{align}
where the constant $\tilde{a}$ depends only on $\calK_\gg$ and
$c_j(\gg)$ is given by \eqref{a24}. 
\end{theorem}

Combining \eqref{1.2.43} and \eqref{1.2.45}, one observes the
remarkable linearity of the theta function with respect to $n$ in
formulas \eqref{1.2.44}, \eqref{1.2.45}. In fact, one can rewrite
\eqref{1.2.45} as
\begin{equation}
b(n)=\Lambda_0
+\sum_{j=1}^{\gg}c_j(\gg)\frac{\partial}{\partial\omega_j}
\ln\bigg(\f{\theta(\underline{\omega}+\ul
A -\ul B n)}{\theta(\underline{\omega}+\ul C -\ul B
n)}\bigg)\bigg|_{\underline{\omega}=0}, \lb{1.3.IM}
\end{equation}
where
\begin{align}
\ul A&= \ul \Xi_{Q_0}-\ua_{Q_0} (P_{\infty+})+\uU_0^{(3)}n_0
+ \ual_{Q_0}(\calD_{\humu (n_0)}),  \lb{1.3.IMA} \\
\ul B&=\uU_0^{(3)}, \lb{1.3.IMB} \\
\ul C&= \ul A + \ul B,  \lb{1.3.IMAA}\\
\Lambda_0&=\frac{1}{2}\sum_{m=0}^{2\gg+1}E_m
-\sum_{j=1}^\gg\lambda_j. \lb{1.3.IML}
\end{align}
Hence, the constants $\Lambda_0\in\bbC$ and $\ul B \in\bbC^\gg$ are
uniquely determined by $\calK_\gg$ (and its homology basis), and
the constant $\ul A\in\bbC^\gg$ is in one-to-one correspondence
with the Dirichlet data
$\humu(n_0)=(\hmu_1(n_0),\dots,\hmu_\gg(n_0)) \in \sym^{\gg}
\calK_\gg$ at the point $n_0$.

\begin{remark} \lb{r2.8}
If one assumes $a$, $b$ in \eqref{1.2.44} and \eqref{1.2.45} to be
quasi-periodic (cf.\ \eqref{3.14a} and \eqref{3.14b}), then there
exists a homology basis $\{\tilde a_j, \tilde b_j\}_{j=1}^\gg$ on
$\calK_\gg$ such that $\wti{\ul B}=\wti{\ul U}^{(3)}_0$ satisfies
the constraint
\begin{equation}
\wti{\ul B}=\wti{\ul U}^{(3)}_0 \in \bbR^\gg. \lb{2.66}
\end{equation}
This is studied in detail in Appendix \ref{sB}.
\end{remark}

\section{The Green's Function of $H$} \label{s3}

In this section we focus on the properties of the Green's function
of $H$ and derive a variety of results to be used in our principal
Section \ref{s4}.

Introducing
\begin{align}
&
\text{G}(P,m,n)=\frac{1}{W\big(\psi(P^*,\cdot,n_0), \psi(P,\cdot,n_0)\big)}
\begin{cases}
\psi(P^*,m,n_0)\psi(P,n,n_0), & \, m\leq n,\\
\psi(P,m,n_0)\psi(P^*,n,n_0), & \, m\geq n,\no
\end{cases}\\
&\hspace*{6.4cm}\quad P\in \calK_\gg \backslash \{
P_{\infty_\pm}\}, \; n,n_0\in\bbZ, \lb{3.4a}
\end{align}
and
\begin{equation}
g(P,n)=G(P,n,n)=\f{\psi(P,n,n_0)\psi(P^*,n,n_0)}{W\big(
\psi(P^*,\cdot,n_0), \psi(P,\cdot,n_0)\big)}, \lb{3.4}
\end{equation}
equations \eqref{1.2.31} and \eqref{1.2.32ab} then imply
\begin{equation}
g(P,n)=-\f{F_\gg(z,n)}{y(P)}, \quad P=(z,y)\in \calK_\gg
\backslash \{ P_{\infty_\pm}\}, \; n\in\bbZ. \lb{3.5}
\end{equation}
Together with $g(P,n)$ we also introduce its two branches
$g_\pm(z,n)$ defined on the upper and lower sheets $\Pi_\pm$ of
$\calK_\gg$ (cf.\ \eqref{a3}, \eqref{a4}, and \eqref{a13})
\begin{equation}
g_\pm (z,n)=\mp \f{F_\gg(z,n)}{R_{2\gg+2}(z)^{1/2}}, \quad
z\in\Pi, \; n\in\bbZ \lb{3.6}
\end{equation}
with $\Pi=\bbC\backslash\calC$ the cut plane introduced in
\eqref{a4}.

For convenience we shall focus on $g_-$ whenever possible and use
the simplified notation
\begin{equation}
g(z,n)=g_-(z,n), \quad z\in\Pi, \; n\in\bbZ \lb{3.7}
\end{equation}
from now on.

Next we briefly review a few properties of quasi-periodic and
almost-periodic discrete functions.

We denote by $\QP(\bbZ)$ and $\AP(\bbZ)$ the sets of quasi-periodic
and almost periodic sequences on $\bbZ$, respectively.

In particular, a sequence $f$ is called quasi-periodic with
fundamental periods $(\Omega_1,\dots,\Omega_N) \in (0,\infty)^N$
if the frequencies $2\pi/\Omega_1,\dots,2\pi/\Omega_N$ are
linearly independent over $\bbQ$ and if there exists a continuous
function $F\in C(\bbR^N)$, periodic of period $1$ in each of its
arguments,
\begin{equation}
F(x_1,\dots,x_j+1,\dots,x_N)=F(x_1,\dots,x_N), \quad x_j\in\bbR,
\;  j=1,\dots,N, \lb{3.14a}
\end{equation}
such that
\begin{equation}
f(n)=F(\Omega_1^{-1}n,\dots,\Omega_N^{-1}n), \quad n\in\bbZ.
\lb{3.14b}
\end{equation}
Any quasi-periodic sequence on $\bbZ$ is almost periodic on
$\bbZ$. Moreover, a sequence $f=\{f(k)\}_{k\in \bbZ}$ is almost periodic
on $\bbZ$ if and only if there exists a Bohr almost periodic
function $g$ on $\bbR$ such that $f(k)=g(k)$ for all $k\in \bbZ$
(see, e.g., \cite[p.\ 47]{Co89}).

For any almost periodic sequence $f=\{f(k)\}_{k\in \bbZ}$, the
mean value $\langle f \rangle$ of $f$, defined by
\begin{equation}
\langle f\rangle =\lim_{N\to\infty}\f{1}{2N+1}
\sum_{k=n_0-N}^{n_0+N} f(k), \lb{3.8}
\end{equation}
exists and is independent of $n_0\in\bbZ$. Moreover, we recall the
following facts for almost periodic sequences that can be deduced
from corresponding properties of Bohr almost periodic functions,
see, for instance, \cite[Ch.\ I]{Be54}, \cite[Sects.\
39--92]{Bo47}, \cite[Ch.\ I]{Co89}, \cite[Chs.\ 1,3,6]{Fi74},
\cite{JM82}, \cite[Chs.\ 1,2,6]{LZ82}, and \cite{Sc65}.

\begin{theorem} \lb{t3.1}
Assume $f,g \in \AP(\bbZ)$ and $n_0,n\in\bbZ$. Then the following
assertions
hold: \\
$(i)$ $f\in \ell^\infty(\bbZ)$.\\
$(ii)$ $\ol f$, $cf$, $c\in\bbC$, $f(\cdot+n)$, $f(n\,\cdot)$,
$n\in\bbZ$,
$|f|^\alpha$, $\alpha\geq 0$ are all in $\AP(\bbZ)$. \\
$(iii)$ $f+g, fg\in \AP(\bbZ)$. \\
$(iv)$ $1/g\in \AP(\bbZ)$ if and only if $1/g\in \ell^\infty(\bbZ)$. \\
$(v)$ Let $G$ be uniformly continuous on $\calM\subseteq \bbR$ and
$f(n)\in\calM$ for all $n\in\bbZ$. Then \\ \hspace*{.5cm} $G(f)\in
\AP(\bbZ)$. \\
$(vi)$ Let $\langle f\rangle=0$, then $\sum_{k=n_0}^n \,
f(k)\underset{|n|\to\infty}{=}\oh(|n|)$.\\
$(vii)$ Let $F(n)=\sum_{k=n_0}^n \, f(k)$. Then
$F\in \AP(\bbZ)$ if and only if $F\in \ell^\infty(\bbZ)$. \\
$(viii)$ If $0\leq f\in \AP(\bbZ)$, $f\not\equiv 0$, then $\langle
f\rangle>0$. \\
$(ix)$ If $1/f\in \ell^\infty(\bbZ)$ and $f=|f|\exp(i\varphi)$,
then $|f|\in \AP(\bbZ)$ and $\varphi$ is of the type \\
\hspace*{.6cm}
$\varphi(n)=cn+\psi(n)$, where $c\in \bbR$ and $\psi\in \AP(\bbZ)$
$($and real-valued\,$)$. \\
$(x)$ If $F(n)=\exp\Big(\sum_{k=n_0}^{n} \, f(k)\Big)$, then $F\in
\AP(\bbZ)$ if and only if $f(n)=i\beta + \psi(n)$, \\
\hspace*{.45cm} where $\beta\in\bbR$, $\psi\in \AP(\bbZ)$, and $\Psi\in
\ell^\infty(\bbZ)$, where $\Psi(n)=\sum_{k=n_0}^{n} \, \psi(k)$.
\end{theorem}

For the rest of this paper it will be convenient to introduce the
following hypothesis:

\begin{hypothesis} \lb{h3.2}
Assume the affine part of $\calK_\gg$ to be nonsingular. Moreover,
suppose that $a,\,b \in \QP(\bbZ)$ satisfy the $\gg$th stationary
Toda equation \eqref{1.2.9} on $\bbZ$.
\end{hypothesis}

Next, we note the following result.

\begin{lemma} \lb{l3.3}
Assume Hypothesis \ref{h3.2}. Then all $z$-derivatives of
$F_{\gg}(z,\cdot)$ and \break $G_{\gg+1}(z,\cdot)$, $z\in\bbC$,
and $g(z,\cdot)$, $z\in\Pi$, are quasi-periodic. Moreover, taking
limits to points on $\calC$, the last result extends to either
side of cuts in the set $\calC\backslash\{E_m\}_{m=0}^{2\gg+1}$
$($cf.\ \eqref{a3}$)$ by continuity with respect to $z$.
\end{lemma}
\begin{proof}
Since $f_\ell$ and $g_\ell$ are polynomials with respect to $a$
and $b$, $f_\ell$ and $g_\ell$, $\ell\in\bbN$, are quasi-periodic
by Theorem \ref{t3.1}. The corresponding assertion for
$F_\gg(z,\cdot)$ then follows from \eqref{1.2.11a} and that for
$g(z,\cdot)$ follows from \eqref{3.6}.
\end{proof}

In the following we represent $G_{\gg+1}(z,n)+G^{+}_{\gg+1}(z,n)$ as
\begin{equation}
G_{\gg+1}(z,n)+G^{+}_{\gg+1}(z,n)=-2\prod_{k=1}^{\gg+1}[z-\nu_k(n)],\quad
z \in \bbC,\; n\in\bbZ, \lb{3.18aa}
\end{equation}
and note that the roots $\nu_{k}$ are bounded,
\begin{equation}
\|\nu_{k}\|_{\infty}\leq\wti C,\quad  k=1,...,\gg+1 \lb{3.4d}
\end{equation}
for some constant $\wti C>0$,  
since the coefficients of $G_{\gg+1}(z,n)$ are defined in terms of
bounded coefficients $a$ and $b$ by\eqref{1.2.4c}. For future
purposes we introduce the set
\begin{align}
\Pi_C&= \Pi  \Big\backslash \Big(\{z\in\bbC\,|\, |z|\leq C+1\} \cup \no \\
& \Big\{z\in\bbC\,|\, \min_{m=0,\dots,2\gg+1}[\Re(E_m)]-1 \leq \Re(z)\leq
\max_{m=0,\dots,2\gg+1}[\Re(E_m)]+1, \no \\
& \quad \min_{m=0,\dots,2\gg+1}[\Im(E_m)]-1\leq \Im(z)\leq
\max_{m=0,\dots,2\gg+1}[\Im(E_m)]+1\Big\}\Big), \lb{3.9}
\end{align}
where $C=\max\{C_{\mu},\|b\|_{\infty}, \wti C\}$ and $C_{\mu}$ is the
constant in \eqref{1.3C}. Without loss of generality, we may also 
assume that $\Pi_C$ contains no cuts, that is,
\begin{equation}
\Pi_C\cap \calC=\emptyset. \lb{3.10}
\end{equation}

Next, we derive a fundamental equation for the mean value of the
diagonal Green's function $g(z,\cdot)$ that will allow us to analyze the
spectrum of the Jacobi operator $H$. First, we note that by
\eqref{1.2.28}, \eqref{1.2.29}, \eqref{1.2.32ab}, and \eqref{3.4a}
one obtains

\begin{align}
& -\frac{\text{G}(P,n,n+1)}{\text{G}(P^{*},n,n+1)}=
\frac{G_{\gg+1}(z,n)-y}{G_{\gg+1}(z,n)+y}, \quad P=(z,y)\in
\calK_{\gg}, \; n\in\bbZ .
     \lb{3.11}
\end{align}
Differentiating the logarithm of the expression on the right-hand
side of \eqref{3.11} with respect to $z$ and using \eqref{1.2.17},
one infers
\begin{align}
& \frac{1}{2}\frac{d}{dz}
\text{ln}\bigg(\frac{G_{\gg+1}(z,n)-y}{G_{\gg+1}(z,n)+y}\bigg)=
\frac{\f{R_{2p+2}^{\bullet}(z)}{2y}G_{\gg+1}(z,n)-y{G}^{\bullet}_{\gg+1}(z,n)}{-4a(n)^2F_\gg(z,n)F_\gg^+(z,n)},\;
z\in \Pi_C. \lb{3.12}
\end{align}
Here $\bullet$ abbreviates $d/dz$. 
We note that the left-hand side of \eqref{3.12} is well-defined since by
\eqref{1.2.17}, \eqref{1.2.21}, and \eqref{3.9}, 
\begin{align}
&[G_{p+1}(z,n)-y][G_{p+1}(z,n)+y] = G_{p+1}(z,n)^2-R_{2p+2}(z) \no\\
& \quad =4a(n)^2F_p(z,n) F_p^+(z,n) \no \\
& \quad =4a(n)^2 \prod_{j=1}^p [z-\mu_j(n)][z-\mu_j(n+1)] \neq 0, \quad 
z\in \Pi_C. 
\end{align}
Adding and subtracting $g(z,n)$ on the right-hand side of \eqref{3.12}
yields
\begin{align}
\frac{1}{2}\frac{d}{dz}
\text{ln}\bigg(\frac{G_{\gg+1}(z,n)-y}{G_{\gg+1}(z,n)+y}\bigg) 
=g(z,n)+\frac{K(z,n)}{y}, \; z\in\Pi_C, \lb{3.13}
\end{align}
where
\begin{equation}
K(z,n)= \f{1}{2}G_{\gg+1}(z,n)
\bigg(\frac{{F}^{\bullet}_\gg(z,n)}{F_\gg(z,n)}+\frac{({F}^+_\gg)^{\bullet}(z,n)}{F_\gg^{+}(z,n)}\bigg)
-{G}^{\bullet}_{\gg+1}(z,n)-F_\gg(z,n). \lb{3.14}
\end{equation}

Next we prove that the mean value of $K(z,\cdot)$ equals zero.

\begin{lemma} \lb{l3.4}
Assume Hypothesis \ref{h3.2}. Then
\begin{equation}
\langle K(z,\cdot) \rangle=0, \quad z\in \Pi_C.
\lb{3.15}
\end{equation}
\end{lemma}
\begin{proof}
Let $z\in\Pi_C$. Using \eqref{1.2.16a} we rewrite \eqref{3.14} as
\begin{align}
K(z,n)=&\f{1}{2}G_{\gg+1}(z,n)\bigg[ \f{d}{dz}\ln \big(
G_{\gg+1}(z,n)+G_{\gg+1}^-(z,n)\big)  \no \\
& \hspace*{2.1cm}+\f{d}{dz}\ln \big(
G_{\gg+1}^+(z,n)+G_{\gg+1}(z,n)\big)
\bigg] \no\\
&-\f{d}{dz}G_{\gg+1}(z,n)+\f{1}{2}\bigg(
\f{G_{\gg+1}^-(z,n)}{z-b(n)}-\f{G_{\gg+1}(z,n)}{z-b^+(n)}\bigg) \no\\
=&\f{1}{2}G_{\gg+1}(z,n)\bigg[
\f{G^{\bullet}_{\gg+1}(z,n)
+(G_{\gg+1}^-)^{\bullet}(z,n)}{G_{\gg+1}(z,n)+G_{\gg+1}^-(z,n)}
\no\\
&
\hspace*{2.1cm}+\f{(G_{\gg+1}^+)^{\bullet}(z,n)
+G^{\bullet}_{\gg+1}(z,n)}{G_{\gg+1}^+(z,n)+G_{\gg+1}(z,n)}
\bigg] \no \\
&-G^{\bullet}_{\gg+1}(z,n) +\f{1}{2}\bigg(
\f{G_{\gg+1}^-(z,n)}{z-b(n)}-\f{G_{\gg+1}(z,n)}{z-b^+(n)}\bigg)\no\\
=&\f{1}{2}\bigg[\f{(G^+_{\gg+1})^{\bullet}(z,n)G_{\gg+1}(z,n)
-G_{\gg+1}^{\bullet}(z,n)G^+_{\gg+1}(z,n)}{G^+_{\gg+1}(z,n)+G_{\gg+1}(z,n)}
\no\\
&\hspace*{.43cm}-\f{G^{\bullet}_{\gg+1}(z,n)G_{\gg+1}^-(z,n)
-(G_{\gg+1}^-)^{\bullet}(z,n)G_{\gg+1}(z,n)}{G_{\gg+1}(z,n)
+G_{\gg+1}^-(z,n)}\bigg]\no\\
&+\f{1}{2}\bigg(
\f{G_{\gg+1}^-(z,n)}{z-b(n)}-\f{G_{\gg+1}(z,n)}{z-b^+(n)}\bigg),\quad
z\in\Pi_C. \lb{3.4f}
\end{align}
Since $K(z,\cdot)$ is a sum of two difference expressions and
$G_{\gg+1}(z,\cdot)$ and $G^{\bullet}_{\gg+1}(z,\cdot)$ are
bounded for fixed $z\in \Pi_C$, one obtains
\begin{equation}
\langle K(z,\cdot) \rangle=0,\, \quad z\in \Pi_C. \lb{3.22a}
\end{equation}
\end{proof}

Using \eqref{3.13} and Lemma \ref{l3.4}, one derives the following
result that will subsequently play a crucial role in this paper.

\begin{lemma} \lb{l3.5}
Assume Hypothesis \ref{h3.2} and let $z, z_0\in\Pi$. Then
\begin{equation}
\bigg<
\ln\bigg(\frac{G_{\gg+1}(z,\cdot)-y}{G_{\gg+1}(z,\cdot)+y}\bigg)\bigg>=2\int_{z_0}^z
dz' \langle g(z',\cdot)\rangle+\bigg<
\ln\bigg(\frac{G_{\gg+1}(z_0,\cdot)-y}{G_{\gg+1}(z_0,\cdot)+y}\bigg)\bigg>,
\lb{3.23}
\end{equation}
where the path connecting $z_0$ and $z$ is assumed to lie in the
cut plane $\Pi$. Moreover, by taking limits to points on $\calC$
in \eqref{3.23}, the result \eqref{3.23} extends to either side of
the cuts in the set $\calC$ by continuity with respect to $z$.
\end{lemma}
\begin{proof}
Let $z, z_0\in\Pi_C$. Integrating equation \eqref{3.13} from $z_0$
to $z$ along a smooth path in $\Pi_C$ yields
\begin{align}
\ln\bigg(\frac{G_{\gg+1}(z,n)-y}{G_{\gg+1}(z,n)+y}\bigg) -
\ln\bigg(\frac{G_{\gg+1}(z_0,\cdot)-y}{G_{\gg+1}(z_0,\cdot)+y}\bigg)
&= 2\int_{z_0}^z dz'\, g(z',n) + \no
\\
& \quad +2\int_{z_0}^z dz'\, \frac{K(z',n)}{y}. \lb{3.18}
\end{align}
By Lemma \ref{l3.3}, $K(z,\cdot)$ is quasi-periodic. Consequently,
also
\begin{equation}
\int_{z_0}^z dz'\, \frac{K(z',\cdot)}{y}, \quad z\in\Pi_C,
\end{equation}
is a family of uniformly almost periodic functions for $z$ varying
in compact subsets of $\Pi_C$ as discussed in \cite[Sect.\ 2.7]{Fi74}.  
By Lemma \ref{l3.4} one thus obtains
\begin{equation}
\bigg\langle \bigg[\int_{z_0}^z dz'\,
\frac{K(z',\cdot)}{y}\bigg]\bigg\rangle =0. \lb{3.21}
\end{equation}
Hence, taking mean values in \eqref{3.18} (taking into account
\eqref{3.21}), proves \eqref{3.23} for $z\in\Pi_C$. Since
$f_\ell$, $\ell\in\bbN_0$, are quasi-periodic by Lemma \ref{l3.3}
(we recall that $f_0=1$), \eqref{1.2.11a} and \eqref{3.6} yield
\begin{equation}
\int_{z_0}^z dz'\, \langle g(z',\cdot)\rangle =
\sum_{\ell=0}^\gg\langle f_{\gg-\ell} \rangle \int_{z_0}^z dz'\,
\f{{(z')}^\ell}{R_{2\gg+2}(z')^{1/2}}.  \lb{3.22}
\end{equation}
Thus, $\int_{z_0}^z dz'\, \langle g(z',\cdot)\rangle$ has an
analytic continuation with respect to $z$ to all of $\Pi$ and
consequently, \eqref{3.23} for $z\in\Pi_C$ extends by analytic
continuation to $z\in\Pi$. By continuity this extends to either
side of the cuts in $\calC$. Interchanging the role of $z$ and
$z_0$, analytic continuation with respect to $z_0$ then yields
\eqref{3.23} for $z,z_0\in\Pi$.
\end{proof}

\begin{remark} \lb{r3.6}
For $z\in\Pi_C$, the sequence
$\text{ln}\big(\frac{G_{\gg+1}(z,\cdot)-y}{G_{\gg+1}(z,\cdot)+y}\big)$
is quasi-periodic and hence
$\big<\text{ln}\big(\frac{G_{\gg+1}(z,\cdot)-y}{G_{\gg+1}(z,\cdot)+y}\big)\big>$
is well-defined. But if one analytically continues \break
$\text{ln}\big(\frac{G_{\gg+1}(z,n)-y}{G_{\gg+1}(z,n)+y}\big)$
with respect to $z$, then $(G_{\gg+1}(z,n)-y)$ and
$(G_{\gg+1}(z,n)+y)$ may acquire zeros for some $n\in\bbZ$ and
hence
$\text{ln}\big(\frac{G_{\gg+1}(z,n)-y}{G_{\gg+1}(z,n)+y}\big)\notin
\QP(\bbZ)$. Nevertheless, as shown by the right-hand side of
\eqref{3.23},
$\big<\text{ln}\big(\frac{G_{\gg+1}(z,\cdot)-y}{G_{\gg+1}(z,\cdot)+y}\big)\big>$
admits an analytic continuation in $z$ from $\Pi_C$ to all of
$\Pi$, and from now on, $\big\langle
\text{ln}\big(\frac{G_{\gg+1}(z,\cdot)-y}{G_{\gg+1}(z,\cdot)+y}\big)\big\rangle$,
$z\in\Pi$, always denotes that analytic continuation (cf.\ also
\eqref{3.25}).
\end{remark}

Next, we will invoke the Baker-Akhiezer function $\psi(P,n,n_0)$
and analyze the expression
$\big<\text{ln}\big(\frac{G_{\gg+1}(z,\cdot)-y}{G_{\gg+1}(z,\cdot)+y}\big)\big>$
in more detail.

\begin{theorem} \lb{t3.7}
Assume Hypothesis \ref{h3.2}, let $P=(z,y)\in \Pi_\pm$, and
$n,n_0\in\bbZ$. Moreover, select a homology basis $\{\ti a_j, \ti
b_j\}_{j=1}^\gg$ on $\calK_\gg$ such that $\wti{\ul B}=\wti{\ul
U}^{(3)}_0$, with $\wti{\ul U}^{(3)}_0$ the vector of $\ti
b$-periods of the normalized differential of the third kind, $\wti
\ome_{P_{\infty_+},P_{\infty_-}}^{(3)}$, satisfies the constraint
\begin{equation}
\wti{\ul B}= \wti{\ul U}^{(3)}_0 \in \bbR^\gg \lb{3.24}
\end{equation}
$($cf.\ Appendix \ref{sB}$)$. Then,
\begin{equation}
\Re\bigg(\bigg\langle
\ln\bigg(\frac{G_{\gg+1}(z,\cdot)-y}{G_{\gg+1}(z,\cdot)+y}\bigg)
\bigg\rangle\bigg) = 2\Re\bigg(\int_{Q_0}^P \wti
\ome_{P_{\infty_+},{\infty_-}}^{(3)}\bigg).  \lb{3.25}
\end{equation}
\end{theorem}
\begin{proof}
Using \eqref{1.2.22a}, \eqref{1.2.25a} and \eqref{1.2.25b} one
obtains the following representation of the Baker-Akhiezer
function $\psi(P,n,n_0)$ for $n>n_0$, $n,n_0\in\bbZ$,
$P\in\calK_{\gg}$,
\begin{align}
     \psi(P,n,n_0)&=\prod_{m=n_0}^{n-1} \phi(P,m)=
\bigg[ \prod_{m=n_0}^{n-1}
\frac{y-G_{\gg+1}(z,m)}{-y-G_{\gg+1}(z,m)}\frac{F_\gg(z,m+1)}{F_\gg(z,m)}\bigg]^{1/2}
\no \\
&
=\bigg(\frac{F_\gg(z,n)}{F_\gg(z,n_0)}\bigg)^{1/2}\bigg[\prod_{m=n_0}^{n-1}
\frac{G_{\gg+1}(z,m)-y}{G_{\gg+1}(z,m)+y}\bigg]^{1/2} \no \\
& =\bigg(\frac{F_\gg(z,n)}{F_\gg(z,n_0)}\bigg)^{1/2} \no
\\
&\quad\times\text{exp}\bigg(\frac{1}{2}\sum_{m=n_0}^{n-1}
\bigg[\text{ln}\bigg(
\frac{G_{\gg+1}(z,m)-y}{G_{\gg+1}(z,m)+y}\bigg) -\bigg<
\text{ln}\bigg(\frac{G_{\gg+1}(z,\cdot)-y}{G_{\gg+1}(z,\cdot)+y}\bigg)\bigg>\bigg]
\bigg)
\no \\
& \quad\times \text{exp}\bigg(\frac{1}{2}(n-n_0)\bigg<
\text{ln}\bigg(\frac{G_{\gg+1}(z,\cdot)-y}{G_{\gg+1}(z,\cdot)+y}\bigg)\bigg>\bigg),
\lb{3.26}\\
& \qquad\qquad\, P=(z,y)\in\Pi_\pm, \; z\in\Pi_C, \; n,n_0\in\bbZ.
\no
\end{align}
A similar representation can be written for $\psi(P,n,n_0)$ if
$n<n_0$, $n,n_0\in\bbZ$, $P\in\calK_{\gg}$.

Since $\Big[\text{ln}\Big(
\frac{G_{\gg+1}(z,m)-y}{G_{\gg+1}(z,m)+y}\Big) -\Big<
\text{ln}\Big(\frac{G_{\gg+1}(z,\cdot)-y}{G_{\gg+1}(z,\cdot)+y}\Big)\Big>\Big]$
has mean zero,
\begin{align}
&\bigg(\frac{1}{2}\sum_{m=n_0}^{n-1} \bigg[\text{ln}\bigg(
\frac{G_{\gg+1}(z,m)-y}{G_{\gg+1}(z,m)+y}\bigg) -\bigg<
\text{ln}\bigg(\frac{G_{\gg+1}(z,\cdot)-y}{G_{\gg+1}(z,\cdot)+y}\bigg)\bigg>\bigg]
\bigg)\underset{|n|\to\infty}{=} \oh(|n|), \no\\
&\hspace*{9.7cm} z\in\Pi_C, \lb{3.27}
\end{align}
by Theorem \ref{t3.1}\,$(vi)$. In addition, the factor $F_\gg
(z,n)/F_\gg (z,n_0)$ in \eqref{3.26} is quasi-periodic and hence
bounded on $\bbZ$.

On the other hand, \eqref{1.2.41} yields
\begin{align}
\psi(P,n,n_0)&=
C(n,n_0)\frac{\theta(\underline{z}(P,\humu(n)))}{\theta(\underline{z}(P,\humu(n_0)))}\exp
\bigg((n-n_0) \int_{Q_0}^P
\omega_{P_{\infty_+},P_{\infty_-}}^{(3)}\bigg) \no \\
&= \Theta(P,n,n_0)\exp \bigg((n-n_0)\int_{Q_0}^P
\wti\ome_{P_{\infty_+},P_{\infty_-}}^{(3)}\bigg), \lb{3.28}
\\
& \hspace*{1.85 cm}
P\in\calK_\gg\backslash
\big(\{P_{\infty_\pm}\}\cup\{\hmu_j(n_0)\}_{j=1}^\gg\big).
\no
\end{align}
Taking into account \eqref{1.2.37}, \eqref{1.2.43}, \eqref{2.66},
\eqref{a31}, and the fact that by \eqref{1.3C} no $\hmu_j(n)$ can
reach $P_{\infty_\pm}$ as $n$ varies in $\bbZ$, one concludes that
\begin{equation}
\Theta(P,\cdot,n_0) \in \ell^\infty(\bbZ), \quad
P\in\calK_\gg\backslash \{\hmu_j(n_0)\}_{j=1}^\gg. \lb{3.29}
\end{equation}
A comparison of \eqref{3.26} and \eqref{3.28} then shows that the
$\oh(|n|)$-term in \eqref{3.27} must actually be bounded on $\bbZ$
and hence the left-hand side of \eqref{3.27} is almost periodic
(in fact, quasi-periodic). In addition, the term
\begin{equation}
\text{exp}\bigg(\frac{1}{2}\sum_{m=n_0}^{n-1}
\bigg[\text{ln}\bigg(
\frac{G_{\gg+1}(z,m)-y}{G_{\gg+1}(z,m)+y}\bigg) -\bigg<
\text{ln}\bigg(\frac{G_{\gg+1}(z,\cdot)-y}{G_{\gg+1}(z,\cdot)+y}\bigg)\bigg>\bigg]
\bigg), \quad z\in\Pi_C, \lb{3.30}
\end{equation}
is then almost periodic (in fact, quasi-periodic) by Theorem
\ref{t3.1}\,$(x)$.\ A further comparison of \eqref{3.26} and
\eqref{3.28} then yields \eqref{3.25} for $z\in\Pi_C$. Analytic
continuation with respect to $z$ then implies \eqref{3.25} for
$z\in\Pi$. By continuity  with respect to $z$, taking boundary
values to either side  of the cuts in the set $\calC$, this then
extends to $z\in\calC$ (cf.\ \eqref{a3}, \eqref{a4}) and hence
proves \eqref{3.25} for $P=(z,y)\in \calK_\gg\backslash \{
P_{\infty_\pm}\}$.
\end{proof}

\section{Spectra of Jacobi Operators with
Quasi-Periodic Algebro-Geometric Coefficients} \lb{s4}

In this section we establish the connection between the
algebro-geometric formalism of Section \ref{s2} and the spectral
theoretic description of Jacobi operators $H$ in $\ell^2(\bbZ)$ with
quasi-periodic algebro-geometric coefficients. In particular, we
introduce the conditional stability set of $H$ and prove our
principal result, the characterization of the spectrum of $H$.
Finally, we provide a qualitative description of the spectrum of
$H$ in terms of analytic spectral arcs.

Suppose that $a$, $b\in \ell^{\infty}(\bbZ)\cap \QP(\bbZ)$ satisfy
the $\gg$th stationary Toda equation \eqref{1.2.9} on $\bbZ$. The
corresponding Jacobi operator $H$ in $\ell^2(\bbZ)$ is then
defined by
\begin{equation}
H=aS^{+}+a^{-}S^{-}+b,\quad \text{dom}(H)=\ell^2(\bbZ). \lb{4.1}
\end{equation}
Thus, $H$ is a bounded operator on $\ell^2(\bbZ)$ (it is
self-adjoint if and only if $a$ and $b$ are real-valued).

Before we turn to the spectrum of $H$ in the general
non-self-adjoint case, we briefly mention the following result on
the spectrum of $H$ in the self-adjoint case with quasi-periodic
(or almost periodic) real-valued coefficients $a$ and $b$. We
denote by $\sigma(A)$, $\sigma_{\rm{e}}(A)$, and $\sigma_{\rm
d}(A)$ the spectrum, essential spectrum, and discrete spectrum of
a self-adjoint operator $A$ in a complex Hilbert space,
respectively.

\begin{theorem} [{\rm See, e.g., \cite{Si82} in the continuous
context}] \lb{t4.1}
Let $a,b\in \QP(\bbZ)$ be real-valued. Define the
self-adjoint Jacobi operator $H$ in $\ell^2(\bbZ)$ as in
\eqref{4.1}. Then,
\begin{align}
&\sigma(H)=\sigma_{\rm{e}}(H)   \\
&\quad\subseteq [-2\sup_{n\in\bbZ}\big(|\Re
(a(n))|\big)+\inf_{n\in\bbZ}\big(\Re(b(n))\big),2\sup_{n\in\bbZ}\big(|\Re
(a(n))|+\sup_{n\in\bbZ}\Re(b(n))\big)], \no\\
&\sigma_{\rm d}(H)=\emptyset. 
\end{align}
Moreover, $\sigma(H)$ contains no isolated points, that is,
$\sigma(H)$ is a perfect set.
\end{theorem}

In the special periodic case where $a,\, b$ are real-valued, the
spectrum  of $H$ is purely absolutely continuous and a finite
union of some compact intervals (see, e.g., \cite{DT176},
\cite{DT276}, \cite{Fl175}, \cite{Te00}--\cite{Mo79},).

Next, we turn to the analysis of the generally non-self-adjoint
operator $H$ in \eqref{4.1}. Assuming Hypothesis \ref{h3.2} we 
introduce the set $\Sigma\subset\bbC$ by
\begin{equation}
\Sigma=\bigg\{\lambda\in\bbC\,\bigg|\, \Re\bigg(\bigg\langle
\text{ln}\bigg(\frac{G_{\gg+1}(\lambda,\cdot)-y}{G_{\gg+1}(\lambda,\cdot)+y}\bigg)\bigg\rangle\bigg)=0\bigg\}.
\lb{4.2}
\end{equation}
Below we will show that $\Sigma$ plays the role of the conditional
stability set of $H$, familiar from the spectral theory of
one-dimensional periodic differential and difference operators.

\begin{lemma} \lb{4.l}
Assume Hypothesis \ref{h3.2}. Then $\Sigma$ coincides with the
conditional stability set of $H$, that is,
\begin{align}
\Sigma&=\{\lambda\in\bbC\,|\, \text{there exists at least one
bounded solution} \no \\
& \hspace*{3.75cm} \text{$0\neq\psi\in \ell^\infty(\bbZ)$ of
$H\psi=\lambda\psi$}\}.  \lb{4.3}
\end{align}
\end{lemma}
\begin{proof}
By \eqref{1.2.41} and \eqref{1.2.42},
\begin{align}
\psi(P,n,n_0)&=
C(n,n_0)\frac{\theta(\underline{z}(P,\humu(n)))}{\theta(\underline{z}(P,\humu(n_0)))}\exp
\bigg((n-n_0) \int_{E_0}^z
\omega_{P_{\infty_+},P_{\infty_-}}^{(3)}\bigg), \lb{4.4} \\
& \hspace*{6.25cm} P=(z,y)\in\Pi_\pm, \no
\end{align}
is a solution of $H\psi=z\psi$ which is bounded on $\bbZ$ if and
only if the exponential function in \eqref{4.4} is bounded on
$\bbZ$. By \eqref{3.25}, the latter holds if and only if
\begin{equation}
\Re\bigg(\bigg\langle
\text{ln}\bigg(\frac{G_{\gg+1}(z,\cdot)-y}{G_{\gg+1}(z,\cdot)+y}\bigg)\bigg\rangle\bigg)=0.
\lb{4.5}
\end{equation}
\end{proof}

\begin{remark} \lb{r3.12}
At first sight our {\it a priori} choice of cuts $\calC$ for
$R_{2\gg+2}(\cdot)^{1/2}$, as described in Appendix \ref{sA},
might seem unnatural as they completely ignore the actual spectrum
of $H$. However, the spectrum of $H$ is not known from the outset,
and in the case of complex-valued potentials, spectral
arcs of $H$ may actually cross each other (cf. Theorem
\ref{t4.3}\,$(iv)$) which renders them unsuitable for cuts of
$R_{2\gg+2}(\cdot)^{1/2}$.
\end{remark}

Before we state our first principal result on the spectrum of $H$,
we find it convenient to recall a number of basic definitions and
well-known facts in connection with the spectral theory of
non-self-adjoint operators (we refer to \cite[Chs.\ I, III,
IX]{EE89}, \cite[Sects.\ 1, 21--23]{Gl65}, \cite[Sects.\ IV.5.6,
V.3.2]{Ka80}, and \cite[p.\ 178--179]{RS78} for more details). Let
$S$ be a densely defined closed operator in complex separable
Hilbert space $\cH$. Denote by $\cB(\cH)$ the Banach space of all
bounded linear operators on $\cH$ and by $\ker(T)$ and $\ran(T)$
the kernel (null space) and range of a linear operator $T$ in
$\cH$. The resolvent set, $\rho(S)$, spectrum, $\sigma(S)$, point
spectrum (the set of eigenvalues), $\sigma_{\rm p}(S)$, continuous
spectrum, $\sigma_{\rm c}(S)$, residual spectrum, $\sigma_{\rm
r}(S)$, field of regularity, $\pi(S)$, approximate point spectrum,
$\sigma_{\rm ap}(S)$, two kinds of essential spectra, $\sigma_{\rm
e}(S)$, and $\wti\sigma_{\rm e}(S)$, the numerical range of $S$,
$\Theta(S)$, and the sets $\Delta (S)$ and $\wti\Delta (S)$ are
defined as follows:
\begin{align}
\rho(S)&=\{z\in\bbC\,|\, (S-z I)^{-1}\in \cB(\cH)\},
\lb{5.15} \\
\sigma(S)&=\bbC\backslash\rho (S), \lb{5.8} \\
\sigma_{\rm p}(S)&=\{\lambda\in\bbC\,|\, \ker(S-\lambda
I)\neq\{0\} \},
\lb{5.9} \\
\sigma_{\rm c}(S)&=\{\lambda\in\bbC \,|\, \text{$\ker(S-\lambda
I)=\{0\}$ and $\ran(S-\lambda I)$ is dense in $\cH$} \no \\ &
\hspace*{6.4cm} \text{but not equal to
$\cH$}\}, \lb{5.10} \\
\sigma_{\rm r}(S)&=\{\lambda\in\bbC\,|\, \text{$\ker(S-\lambda
I)=\{0\}$
and $\ran(S-\lambda I)$ is not dense in $\cH$}\}, \lb{5.11} \\
\pi(S)&=\{z \in\bbC \,|\, \text{there exists $k_z >0$ s.t.
$\| (S- zI)u\|_\cH \ge k_z \| u\|_\cH$} \no \\
& \hspace*{6.15cm} \text{for all $u\in\dom(S)$}\}, \lb{5.14} \\
\sigma_{\rm ap}(S)&=\bbC\backslash\pi(S), \lb{5.21} \\
\Delta(S)&=\{z\in\bbC \,|\, \text{$\dim(\ker(S-zI))<\infty$ and
$\ran(S-zI)$ is closed}\}, \lb{5.16} \\
\sigma_{\rm e}(S)&= \bbC\backslash\Delta (S), \lb{5.22b} \\
\wti\Delta(S)&=\{z\in\bbC \,|\, \text{$\dim(\ker(S-zI))<\infty$ or
$\dim(\ker(S^*-\ol z I))<\infty$}\}, \lb{5.16a} \\
\wti\sigma_{\rm e} (S)&=\bbC\backslash \wti\Delta(S), \lb{5.16b} \\
\Theta(S)&=\{(f,Sf)\in\bbC \,|\, f\in\dom(S), \, \|f\|_{\cH}=1\},
\lb{5.16c}
\end{align}
respectively. One then has
\begin{align}
\sigma (S)&=\sigma_{\rm p}(S)\cup\sigma_{\rm{c}}(S)\cup
\sigma_{\rm r}(S) \quad \text{(disjoint union)} \lb{5.17} \\
&=\sigma_{\rm p}(S)\cup\sigma_{\rm{e}}(S)\cup\sigma_{\rm r}(S),
\lb{5.18} \\
\sigma_{\rm c}(S)&\subseteq\sigma_{\rm e}(S)\backslash
(\sigma_{\rm p}(S)\cup\sigma_{\rm r}(S)), \lb{5.18a} \\
\sigma_{\rm r}(S)&=\sigma_{\rm p}(S^*)^* \backslash\sigma_{\rm
p}(S),
\lb{5.19} \\
\sigma_{\rm ap}(S)&=\{\lambda \in \bbC \, |\, \text{there exists a
sequence $\{ f_n\}_{n\in\bbN}\subset\dom(S)$} \no \\
&\hspace*{.1cm} \text{with $\| f_n \|_\cH=1$, $n\in\bbN$, and
$\lim_{n\to\infty} \|(S-\lambda I)f_n\|_\cH=0$}\}, \lb{5.12} \\
\wti\sigma_{\rm e}(S)&\subseteq \sigma_{\rm
e}(S)\subseteq\sigma_{\rm ap}(S)\subseteq\sigma(S) \, \text{ (all
four sets are closed)}, \lb{5.23}
\\
\rho(S)&\subseteq \pi(S) \subseteq \Delta (S) \subseteq \wti\Delta
(S)
\;\; \text{ (all four sets are open),} \lb{5.24} \\
\wti\sigma_{\rm e}(S) & \subseteq \ol{\Theta(S)}, \quad \Theta(S)
\,
\text{ is convex,} \lb{5.25} \\
\wti\sigma_{\rm e}(S) &=\sigma_{\rm e}(S) \, \text{ if $S=S^*$.}
\lb{5.26}
\end{align}
Here $\sigma^*$ in the context of \eqref{5.19} denotes the complex
conjugate of the set $\sigma\subseteq\bbC$, that is,
\begin{equation}
\sigma^*=\{\ol{\lambda}\in\bbC \,|\, \lambda\in\sigma\}. \lb{5.27}
\end{equation}
We note that there are several other versions of the concept of
the essential spectrum in the non-self-adjoint context (cf.\
\cite[Ch.\ IX]{EE89}) but we will only use the two in
\eqref{5.22b} and in \eqref{5.16b} in this paper.

We start with the following elementary result.

\begin{lemma}  \lb{l3.10a}
Let $H$ be defined as in \eqref{4.1}. Then,
\begin{equation}
\sigma_{\rm e}(H)=\wti\sigma_{\rm e}(H)\subseteq \ol{\Theta(H)}.
\lb{3.39a}
\end{equation}
\end{lemma}
\begin{proof}
Since $H$ and $H^*$ are second-order difference operators on
$\bbZ$,
\begin{equation}
\dim(\ker(H-z I))\leq 2, \quad \dim(\ker(H^*-\ol z I))\leq 2.
\lb{3.39b}
\end{equation}
Moreover, we note that $S$ closed and densely defined and
$\dim(\ker(S^*- \ol z I))<\infty$ implies that $\ran(S-zI)$ is
closed (cf.\ \cite[Theorem\ I.3.2]{EE89}). Equations
\eqref{5.16}--\eqref{5.16b} and \eqref{5.25} then prove
\eqref{3.39a}.
\end{proof}

\begin{theorem}  \lb{t3.11}
Assume Hypothesis \ref{h3.2}. Then the point spectrum and residual
spectrum of $H$ are empty and hence the spectrum of $H$ is purely
continuous,
\begin{align}
&\sigma_{\rm p}(H)=\sigma_{\rm r}(H)=\emptyset,  \lb{3.40} \\
&\sigma(H)=\sigma_{\rm c}(H)=\sigma_{\rm e}(H)=\sigma_{\rm ap}(H).
\lb{3.41}
\end{align}
\end{theorem}
\begin{proof}
First we prove the absence of the point spectrum of $H$. Suppose
$z\in\Pi\backslash(\Sigma \cup\{\mu_j(n_0)\}_{j=1}^\gg)$. Then
$\psi(P,\cdot,n_0)$ and $\psi(P^*,\cdot,n_0)$ are linearly
independent solutions of $H\psi=z\psi$ which are unbounded at
$+\infty$ or $-\infty$. This argument extends to all
$z\in\Pi\backslash\Sigma$ by multiplying $\psi(P,\cdot,n_0)$ and
$\psi(P^*,\cdot,n_0)$ with an appropriate function of $z$ and
$n_0$ (independent of $n$). It also extends to either side of the
cut $\calC\backslash\Sigma$ by continuity with respect to $z$. On
the other hand, any solution $\psi(z,\cdot)\in \ell^2(\bbZ)$ of
$H\psi=z\psi$, $z\in\bbC$, is necessarily bounded (since any
sequence in $\ell^2(\bbZ)$ is bounded). Thus,
\begin{equation}
\{\bbC\backslash\Sigma\}\cap \sigma_{\rm p}(H)=\emptyset.
\lb{3.67a}
\end{equation}
Hence, it remains to rule out eigenvalues located in $\Sigma$. We
consider a fixed $\lambda\in\Sigma$ and note that by
\eqref{1.2.31}, there exists at least one solution
$\psi_1(\lambda,\cdot) \in \ell^\infty(\bbZ)$ of
$H\psi=\lambda\psi$. Actually, a comparison of \eqref{3.26} and
\eqref{4.2} shows that one can choose $\psi_1(\lambda,\cdot)$ such
that $|\psi_1(\lambda,\cdot)|\in \QP(\bbZ)$ and hence
$\psi_1(\lambda,\cdot)\notin \ell^2(\bbZ)$.

Next, suppose there exists a second solution
$\psi_2(\lambda,\cdot)\in \ell^2(\bbZ)$ of $H\psi=\lambda\psi$
which is linearly independent of $\psi_1(\lambda,\cdot)$. Then one
concludes that the Wronskian of $\psi_1(\lambda,\cdot)$ and
$\psi_2(\lambda,\cdot)$ lies in $\ell^2(\bbZ)$,
\begin{equation}
W(\psi_1(\lambda,\cdot),\psi_2(\lambda,\cdot))\in \ell^2(\bbZ).
\lb{3.68c}
\end{equation}
However, by hypothesis,
$W(\psi_1(\lambda,\cdot),\psi_2(\lambda,\cdot))=c(\lambda)\neq 0$
is a nonzero constant. This contradiction proves that
\begin{equation}
\Sigma\cap \sigma_{\rm p}(H)=\emptyset \lb{3.68d}
\end{equation}
and hence $\sigma_{\rm p}(H)=\emptyset$.

Next, we note that the same argument yields that $H^*$ also has no
point spectrum,
\begin{equation}
\sigma_{\rm p}(H^*)=\emptyset. \lb{3.43}
\end{equation}
Indeed, if $a,\, b \in \ell^\infty(\bbZ)\cap \QP(\bbZ)$ satisfy the
$\gg$th stationary Toda equation \eqref{1.2.9} on $\bbZ$, then
$\ol a,\, \ol b$ also satisfy one of the $\gg$th stationary Toda
equation \eqref{1.2.9} associated with a hyperelliptic curve of
genus $\gg$ with $\{E_m\}_{m=0}^{2\gg+1}$ replaced by $\{\ol
E_m\}_{m=0}^{2\gg+1}$, etc. Since by general principles (cf.\
\eqref{5.27}),
\begin{equation}
\sigma_{\rm r}(B)\subseteq \sigma_{\rm p}(B^*)^* \lb{3.44}
\end{equation}
for any densely defined closed linear operator $B$ in some complex
separable Hilbert space (see, e.g., \cite[p.\ 71]{Go85}), one
obtains $\sigma_{\rm r}(H)=\emptyset$ and hence \eqref{3.40}. This
proves that the spectrum of $H$ is purely continuous,
$\sigma(H)=\sigma_{\rm c}(H)$. The remaining equalities in
\eqref{3.41} then follow from \eqref{5.18a} and \eqref{5.23}.
\end{proof}

The following result is a fundamental one:

\begin{theorem}  \lb{t4.2}
Assume Hypothesis \ref{h3.2}. Then the spectrum of $H$ coincides
with $\Sigma$ and hence equals the conditional stability set
of $H$,
\begin{align}
\sigma(H) &=\bigg\{\lambda\in\bbC\,\bigg|\, \Re\bigg(\bigg\langle
\ln\bigg(\frac{G_{\gg+1}(\lambda,\cdot)-y}
{G_{\gg+1}(\lambda,\cdot)+y}\bigg)\bigg\rangle\bigg)=0\bigg\}
\lb{4.6} \\
&=\{\lambda\in\bbC\,|\, \text{there exists at least one bounded
solution}  \no \\
& \hspace*{3.7cm} \text{$0\neq\psi\in \ell^\infty(\bbZ)$ of
$H\psi=\lambda\psi$}\}.  \lb{4.7}
\end{align}
In particular,
\begin{equation}
\{E_m\}_{m=0}^{2\gg+1}\subset\sigma(H), \lb{4.8}
\end{equation}
and $\sigma(H)$ contains no isolated points.
\end{theorem}
\begin{proof}
First we will prove that
\begin{equation}
\sigma(H)\subseteq \Sigma \lb{3.42}
\end{equation}
by adapting a method due to Chisholm and Everitt \cite{CE70} (in
the context of differential operators). For this purpose we
temporarily choose
$z\in\Pi\backslash(\Sigma\cup\{\mu_j(n_0)\}_{j=1}^\gg)$ and
construct the resolvent of $H$ as follows. Introducing the two
branches $\psi_\pm (P,n,n_0)$ of the Baker--Akhiezer function
$\psi(P,n,n_0)$ by
\begin{equation}
\psi_\pm (P,n,n_0)=\psi(P,n,n_0), \quad P=(z,y)\in\Pi_\pm, \;
n,n_0\in\bbZ, \lb{3.47}
\end{equation}
we define
\begin{align}
\hat \psi_+(z,n,n_0)&=\begin{cases} \psi_+(z,n,n_0) & \text{if
$\psi_+(z,\cdot,n_0)\in \ell^2(n_0,\infty)$,} \\
\psi_-(z,n,n_0) & \text{if
$\psi_-(z,\cdot,n_0)\in \ell^2(n_0,\infty)$,} \end{cases} \lb{3.48} \\
\hat \psi_-(z,n,n_0)&=\begin{cases} \psi_-(z,n,n_0) & \text{if
$\psi_-(z,\cdot,n_0)\in \ell^2(-\infty,n_0)$,} \\
\psi_+(z,n,n_0) & \text{if
$\psi_+(z,\cdot,n_0)\in \ell^2(-\infty,n_0)$,} \end{cases} \lb{3.49} \\
& \hspace*{4.1cm} z\in\Pi\backslash\Sigma, \; n,n_0\in\bbZ, \no
\end{align}
and
\begin{align}
G(z,n,n')&=\f{1}{W(\hat\psi_-(z,n,n_0),\hat\psi_+(z,n,n_0))}\begin{cases}
\hat\psi_-(z,n',n_0)\hat\psi_+(z,n,n_0), & n\geq n', \\
\hat\psi_-(z,n,n_0)\hat\psi_+(z,n',n_0), & n\leq n', \end{cases} \no \\
& \hspace*{6.4cm} z\in\Pi\backslash\Sigma, \; n,n_0\in\bbZ.
\lb{3.50}
\end{align}
Due to the homogeneous nature of $G$, \eqref{3.50} extends to all
$z\in\Pi$. Moreover, we extend \eqref{3.48}--\eqref{3.50} to
either side of the cut $\calC$ except at possible points in
$\Sigma$ (i.e., to $\calC\backslash\Sigma$) by continuity with
respect to $z$, taking limits to $\calC\backslash\Sigma$. Next, we
introduce the operator $R(z)$ in $\ell^2(\bbZ)$ defined by
\begin{align}
(R(z)f)(n)=\sum_{n'\in\bbZ} \,G(z,n,n')f(n'), \quad f\in
\ell^\infty_0(\bbZ), \; z\in\Pi, \lb{3.51}
\end{align}
where $\ell^{\infty}_0(\bbZ)$ denotes the linear space of
compactly supported (i.e., finite) complex-valued sequences, and
extend it to $z\in\calC\backslash\Sigma$, as discussed in
connection with $G(\cdot,n,n')$. The explicit form of
$\hat\psi_\pm(z,n,n_0)$, inferred from \eqref{3.28} by restricting
$P$ to $\Pi_\pm$, then yields the estimates
\begin{equation}
|\hat \psi_\pm(z,n,n_0)|\leq C_\pm(z,n_0) e^{\mp \kappa(z)n},
\quad z\in\Pi\backslash\Sigma, \;  n\in\bbZ \lb{3.51a}
\end{equation}
for some constants $C_\pm(z,n_0)>0$, $\kappa(z)>0$,
$z\in\Pi\backslash\Sigma$. One can follow the second part of the
proof of Theorem\ 5.3.2 in \cite{Ea73} line by line and prove that
$R(z)$, $z\in\bbC\backslash\Sigma$, extends from
$\ell^\infty_0(\bbZ)$ to a bounded linear operator defined on all
of $\ell^2(\bbZ)$.
A straightforward computation then proves
\begin{equation}
(H-zI)R(z)f=f, \quad f\in \ell^2(\bbZ), \;
z\in\bbC\backslash\Sigma \lb{3.52}
\end{equation}
and hence also
\begin{equation}
R(z)(H-zI)g=g, \quad g\in \ell^2(\bbZ), \;
z\in\bbC\backslash\Sigma. \lb{3.53}
\end{equation}
Thus, $R(z)=(H-zI)^{-1}$, $z\in\bbC\backslash\Sigma$, and hence
\eqref{3.42} holds.

Next we will prove that
\begin{equation}
\sigma(H)\supseteq \Sigma. \lb{3.54}
\end{equation}
We will adapt a strategy of proof applied by Eastham in the
continuous case of (real-valued) periodic potentials \cite{Ea67}
(reproduced in the proof of Theorem\ 5.3.2 of \cite{Ea73}) to the
(complex-valued) quasi-periodic discrete case at hand. Suppose
$\lambda\in\Sigma$. By the characterization \eqref{4.3} of
$\Sigma$, there exists a bounded solution $\psi(\lambda,\cdot)$ of
$H\psi=\lambda\psi$. A comparison with the Baker-Akhiezer function
\eqref{3.28} then shows that one can assume, without loss of
generality, that
\begin{equation}
|\psi(\lambda,\cdot)|\in \QP(\bbZ). \lb{3.55}
\end{equation}
By Theorem \ref{t3.1}\,$(i)$, one obtains
\begin{equation}
\psi(\lambda,\cdot)\in \ell^\infty(\bbZ). \lb{3.55a}
\end{equation}

Next, we pick $\Omega \in \bbN$ and consider $g(n), \;
n=0,1,\dots,\Omega$, satisfying
\begin{align}
\begin{split}
&g(0)=0, \quad g(\Omega)=1, \\
&0\leq g(n)\leq 1, \quad n=1,\dots,\Omega-1. \lb{3.55b}
\end{split}
\end{align}
Moreover, we introduce the sequence $\{h_k\}_{k\in\bbN}\in
\ell^2(\bbZ)$ by
\begin{equation}
h_k(n)=\begin{cases} 1, & |n|\leq (k-1)\Omega, \\
g(k\Omega-|n|), & (k-1)\Omega\leq |n|\leq k\Omega, \\
0, & |n|\geq k\Omega \end{cases} \lb{3.55c}
\end{equation}
and the sequence $\{f_k(\lambda)\}_{k\in\bbN}\in \ell^2(\bbZ)$ by
\begin{equation}
f_k(\lambda,n)=d_k(\lambda)\psi(\lambda,n)h_k(n), \quad n\in\bbZ,
\; d_k(\lambda)>0, \; k\in\bbN. \lb{3.55d}
\end{equation}
Here $d_k(\lambda)$ is determined by the normalization requirement
\begin{equation}
\|f_k(\lambda)\|_2 =1, \quad k\in\bbN. \lb{3.55e}
\end{equation}
Of course,
\begin{equation}
f_k(\lambda,\cdot) \in \ell^2(\bbZ), \quad k\in\bbN, \lb{3.55ea}
\end{equation}
since $f_k(\lambda,\cdot)$ is finitely supported. Next, we note
that as a consequence of Theorem \ref{t3.1}\,$(viii)$,
\begin{equation}
\sum_{-N}^N  |\psi(\lambda,n)|^2\underset{N\to\infty}{=}
(2N+1)\big\langle |\psi(\lambda,\cdot)|^2 \big\rangle +\oh(N)
\lb{3.55f}
\end{equation}
with
\begin{equation}
\big\langle |\psi(\lambda,\cdot)|^2 \big\rangle >0. \lb{3.55g}
\end{equation}
Thus, one computes
\begin{align}
1&=\|f_k(\lambda)\|^2_2=d_k(\lambda)^2\sum_{n \in \bbZ}
|\psi(\lambda,n)|^2
h_k(n)^2 \no \\
& = d_k(\lambda)^2\sum_{|n|\leq k\Omega} |\psi(\lambda,n)|^2
h_k(n)^2 \geq d_k(\lambda)^2\sum_{|n|\leq (k-1)\Omega}
|\psi(\lambda,n)|^2 \no \\
& \quad \geq d_k(\lambda)^2 \big[\big\langle
|\psi(\lambda,\cdot)|^2 \big\rangle (k-1)\Omega + \oh(k) \big].
\lb{3.55h}
\end{align}
Consequently,
\begin{equation}
d_k(\lambda)\underset{k\to\infty}{=} \Oh\big(k^{-1/2}\big).
\lb{3.55i}
\end{equation}
Next, one computes
\begin{align}
(H-\lambda I)f_k(\lambda,n)=&
d_k(\lambda)\Big[a(n)\psi(\lambda,n)\big[h_k(n+1)-h_k(n)\big]  \no \\
&  +a(n-1)\psi(\lambda,n-1)\big[h_k(n-1)-h_k(n)\big]\Big]
\lb{3.55j}
\end{align}
and hence
\begin{align}
\|(H-\lambda I)f_k\|_2 & \leq
2d_k(\lambda)\|a\|_{\infty}\|\psi(\lambda)(h_k^+ - h_k)\|_2, \quad
k\in\bbN. \lb{3.55k}
\end{align}
Using \eqref{3.55a} and \eqref{3.55c} one estimates
\begin{align}
\|\psi(\lambda)\big[h_k^+ - h_k\big]\|_2^2&= \sum_{(k-1)\Omega\leq
|n|\leq k\Omega} \, |\psi(\lambda,n)|^2 |h_k(n+1)-h_k(n)|^2 \no\\
& \leq 2\|\psi(\lambda)\|_{\infty}^2 \big(\Omega+1\big).
\lb{3.55l}
\end{align}
Thus, combining \eqref{3.55i} and \eqref{3.55k}--\eqref{3.55l} one
infers
\begin{equation}
\lim_{n\to\infty}\|(H-\lambda I)f_k\|_2 =0 \lb{3.55n}
\end{equation}
and hence $\lambda\in\sigma_{\rm ap}(H)=\sigma(H)$ by \eqref{5.12}
and \eqref{3.41}.

Relation \eqref{4.8} follows from \eqref{4.3} and the fact that by
\eqref{1.2.31} there exists a solution $\psi((E_m,0),\cdot,n_0)\in
\ell^\infty(\bbZ)$ of $H\psi=E_m\psi$ for all $m=0,\dots,2\gg+1$.

Finally, $\sigma(H)$ contains no isolated points since those would
necessarily be essential singularities of the resolvent of $H$, as
$H$ has no eigenvalues by \eqref{3.40} (cf.\ \cite[Sect.\
III.6.5]{Ka80}). An analysis of the Green's function of $H$
reveals at most an algebraic singularity at the points
$\{E_m\}_{m=0}^{2\gg+1}$ and hence excludes the possibility of an
essential singularity of $(H-zI)^{-1}$.
\end{proof}

In the special self-adjoint case where $a,\,b$ are real-valued,
the result \eqref{4.6} is equivalent to the vanishing of the
Lyapunov exponent of $H$ which characterizes the (purely
absolutely continous) spectrum of $H$ as discussed by
Carmona and Lacroix \cite[Chs.\ IV, VII]{CL90} (cf.\ also \cite{CK87},
\cite{DS83}, \cite{GK03}, \cite{KK03}).

The explicit formula for $\Sigma$ in \eqref{4.2} permits a
qualitative description of the spectrum of $H$ as follows. We
recall \eqref{3.13} and \eqref{3.22} and write
\begin{equation}
\frac{1}{2}\frac{d}{dz}
\bigg\langle\text{ln}\bigg(\frac{G_{\gg+1}(z,\cdot)-y}
{G_{\gg+1}(z,\cdot)+y}\bigg)\bigg\rangle
=\langle
g(z,\cdot)\rangle=\frac{\prod_{j=1}^\gg(z-\wti\lambda_j)}{\Big(\prod_{j=0}^{2\gg+1}(z-E_m)\Big)^{1/2}},\quad
z\in\Pi, \lb{4.9}
\end{equation}
for some constants
\begin{equation}
\{\wti\lambda_j\}_{j=1}^{\gg}\subset\bbC. \lb{4.10}
\end{equation}
As in similar situations before, \eqref{4.9} extends to either
side of the cuts in $\calC$ by continuity with respect to $z$.

\begin{theorem}  \lb{t4.3}
Assume Hypothesis \ref{h3.2}. Then the spectrum $\sigma(H)$ of $H$
has
the following properties: \\
$(i)$ $\sigma(H)\subset \bbC$ is bounded,
\begin{equation}
\sigma(H)\subset \{z\in\bbC\,|\, \Re(z)\in [M_1,M_2],\, \Im(z)\in
[M_3,M_4]\}, \lb{4.11}
\end{equation}
where
\begin{align}
\begin{split}
& M_1=-2\sup_{n\in\bbZ}[|\Re (a(n))|]+\inf_{n\in\bbZ}[\Re(b(n))],\\
& M_2=2\sup_{n\in\bbZ}[|\Re (a(n))|]+\sup_{n\in\bbZ}[\Re(b(n))], \\
& M_3=-2\sup_{n\in\bbZ}[|\Im (a(n))|]+\inf_{n\in\bbZ}[\Im(b(n))],\\
& M_4=2\sup_{n\in\bbZ}[|\Im (a(n))|]+\sup_{n\in\bbZ}[\Im(b(n))].
\lb{4.12}
\end{split}
\end{align}
$(ii)$ $\sigma(H)$ consists of finitely many simple analytic arcs
$($cf. Remark \ref{r4.2}$)$. These analytic arcs may only end at
the points
$\wti\lambda_1,\dots,\wti\lambda_\gg$, $E_0,\dots,E_{2\gg+1}$. \\
$(iii)$ Each $E_m$, $m=0,\dots,2\gg+1$, is met by at least one of
these arcs. More precisely, a particular $E_{m_0}$ is hit by
precisely $2N_0+1$ analytic arcs, where $N_0\in\{0,\dots,\gg\}$
denotes the number of $\wti\lambda_j$ that coincide with
$E_{m_0}$. Adjacent arcs meet at an angle $2\pi/(2N_0+1)$ at
$E_{m_0}$. $($Thus,
generically, $N_0=0$ and precisely one arc hits $E_{m_0}$.$)$ \\
$(iv)$ Crossings of spectral arcs are permitted. This phenomenon
takes place precisely when for a particular
$j_0\in\{1,\dots,\gg\}$, $\wti\lambda_{j_0}\in\sigma(H)$ such that
\begin{align}
\begin{split}
&\Re\bigg(\bigg\langle
\ln\bigg(\frac{G_{\gg+1}(\wti\lambda_{j_0},\cdot)-y}
{G_{\gg+1}(\wti\lambda_{j_0},\cdot)+y}\bigg)\bigg\rangle\bigg)=0
\lb{4.14} \\
& \quad \text{ for some $j_0\in\{1,\dots,\gg\}$ with
$\wti\lambda_{j_0}\notin \{E_m\}_{m=0}^{2\gg+1}$}.
\end{split}
\end{align}
In this case $2M_0+2$ analytic arcs are converging toward
$\wti\lambda_{j_0}$, where $M_0\in\{1,\dots,\gg\}$ denotes the
number of $\wti\lambda_j$ that coincide with $\wti\lambda_{j_0}$.
Adjacent arcs meet at an angle $\pi/(M_0+1)$ at
$\wti\lambda_{j_0}$. $($Thus, if crossings occur, generically,
$M_0=1$ and two arcs cross at a right angle.$)$ \\
$(v)$ The resolvent set $\bbC\backslash\sigma (H)$ of $H$ is
path-connected.
\end{theorem}
\begin{proof}
Item $(i)$ follows from \eqref{3.39a} and \eqref{3.41} upon
noticing that
\begin{equation}
(f,Hf)=2\sum_{k=-\infty}^{\infty}
a(k)\Re[f(k+1)\overline{f(k)}]+(f,\Re(b)f)+i(f,\Im(b)f), \quad
f\in \ell^{2}(\bbZ). \lb{4.15}
\end{equation}

To prove $(ii)$ we first introduce the meromorphic differential of
the third kind
\begin{align}
&\Omega^{(3)} = \langle g(P,\cdot)\rangle dz= \f{\langle
F_\gg(z,\cdot)\rangle dz}{y}=\f{\prod_{j=1}^\gg \big(z-\wti
\lambda_j\big) dz}{R_{2\gg+2}(z)^{1/2}}, \no\\
&\hspace{4.5cm}  P=(z,y)\in \calK_\gg\backslash\{P_{\infty_\pm}\}
\lb{4.16}
\end{align}
(cf.\ \eqref{4.10}). Then, by Lemma \ref{l3.5},
\begin{equation}
\bigg<
\text{ln}\bigg(\frac{G_{\gg+1}(z,\cdot)-y}{G_{\gg+1}(z,\cdot)+y}\bigg)\bigg>=2\int_{Q_0}^P
\Omega^{(3)}+\bigg<
\text{ln}\bigg(\frac{G_{\gg+1}(z_0,\cdot)-y}{G_{\gg+1}(z_0,\cdot)+y}\bigg)\bigg>,
\lb{4.17}
\end{equation} for some fixed $Q_0=(z_0,y_0)\in
\calK_\gg\backslash\{P_{\infty_\pm}\}$, is holomorphic on
$\calK_\gg\backslash\{P_{\infty_\pm}\}$. By \eqref{4.9},
\eqref{4.10}, the characterization \eqref{4.6} of the spectrum,
\begin{equation}
\sigma(H) = \bigg\{\lambda\in\bbC\,\bigg|\, \Re\bigg(\bigg\langle
\ln\bigg(\frac{G_{\gg+1}(\lambda,\cdot)-y}{G_{\gg+1}(\lambda,\cdot)+y}\bigg)\bigg\rangle\bigg)=0\bigg\},
\lb{4.18}
\end{equation}
and the fact that $\Re\big(\big\langle
\ln\big(\frac{G_{\gg+1}(z,\cdot)-y}{G_{\gg+1}(z,\cdot)+y}\big)\big\rangle\big)$
is a harmonic function on the cut plane $\Pi$, the spectrum
$\sigma(H)$ of $H$ consists of analytic arcs which may only end at
the points $\wti\lambda_1,\dots,\wti\lambda_\gg$,
$E_0,\dots,E_{2\gg+1}$. (Since $\sigma(H)$ is independent of the
chosen set of cuts, if a spectral arc crosses or runs along a part
of one of the cuts in $\calC$, one can slightly deform the
original set of cuts to extend an analytic arc along or across
such an original cut.)

To prove $(iii)$ one first recalls that by Theorem \ref{t4.2} the
spectrum of $H$ contains no isolated points. On the other hand,
since $\{E_m\}_{m=0}^{2\gg+1}\subset\sigma(H)$ by \eqref{4.8}, one
concludes that at least one spectral arc meets each $E_m$,
$m=0,\dots,2\gg+1$. Choosing $Q_0=(E_{m_0},0)$ in \eqref{4.17} one
obtains
\begin{align}
&\bigg\langle
\ln\bigg(\frac{G_{\gg+1}(z,\cdot)-y}{G_{\gg+1}(z,\cdot)+y}\bigg)\bigg\rangle
= 2\int_{E_{m_0}}^z dz' \, \langle g(z',\cdot)\rangle +
\bigg\langle
\ln\bigg(\frac{G_{\gg+1}(E_{m_0},\cdot)-y}
{G_{\gg+1}(E_{m_0},\cdot)+y}\bigg)\bigg\rangle
   \no \\
&\quad= {2}\int_{E_{m_0}}^z dz' \f{\prod_{j=1}^\gg
\big(z'-\wti\lambda_j\big)}{\big(\prod_{m=0}^{2\gg+1} (z'-E_m)
\big)^{1/2}}+ \bigg\langle
\ln\bigg(\frac{G_{\gg+1}(E_{m_0},\cdot)-y}
{G_{\gg+1}(E_{m_0},\cdot)+y}\bigg)\bigg\rangle
\no \\
&\underset{z\to E_{m_0}}{=} \int_{E_{m_0}}^z dz'\,
(z'-E_{m_0})^{N_0-(1/2)}[C+\Oh(z'-E_{m_0})] \no\\
&\quad\quad +\bigg\langle
\ln\bigg(\frac{G_{\gg+1}(E_{m_0},\cdot)-y}
{G_{\gg+1}(E_{m_0},\cdot)+y}\bigg)\bigg\rangle
\no \\
&\underset{z\to E_{m_0}}{=}
\frac{(z-E_{m_0})^{N_0+(1/2)}}{N_0+(1/2)}[C+\Oh(z-E_{m_0})]+
\bigg\langle
\ln\bigg(\frac{G_{\gg+1}(E_{m_0},\cdot)-y}{G_{\gg+1}(E_{m_0},\cdot)
+y}\bigg)\bigg\rangle  \lb{4.19}
\end{align}
for some $C=|C|e^{i\varphi_0}\in\bbC\backslash\{0\}$. Using
\begin{equation}
\Re\bigg(\bigg\langle
\ln\bigg(\frac{G_{\gg+1}(E_{m},\cdot)-y}{G_{\gg+1}(E_{m},\cdot)+y}\bigg)
\bigg\rangle\bigg) =0, \quad m=0,\dots,2\gg+1, \lb{4.20}
\end{equation}
as a consequence of \eqref{4.8}, $\Re\big(\big\langle
\ln\big(\frac{G_{\gg+1}(z,\cdot)-y}{G_{\gg+1}(z,\cdot)+y}\big)\big\rangle\big)=0$
and $z=E_{m_0}+\rho e^{i\varphi}$ imply
\begin{equation}
0\underset{\rho\downarrow 0}{=}
\cos[(N_0+(1/2))\varphi+\varphi_0]\rho^{N_0+(1/2)}[|C|+\Oh(\rho)].
\lb{4.21}
\end{equation}
This proves the assertions made in item $(iii)$.

In order to prove $(iv)$ it suffices to refer to \eqref{4.9} and
     observe that locally $\frac{1}{2}\frac{d}{dz}
\big\langle\text{ln}\big(\frac{G_{\gg+1}(z,\cdot)-y}{G_{\gg+1}(z,\cdot)+y}\big)\big\rangle$
behaves like $C_0(z-\wti\lambda_{j_0})^{M_0}$ for some
$C_0\in\bbC\backslash\{0\}$ in a sufficiently small neighborhood
of $\wti\lambda_{j_0}$.

Finally we will show that all arcs are simple (i.e., do not
self-intersect each other). Assume that the spectrum of $H$
contains a simple closed loop $\gamma$, $\gamma\subset\sigma(H)$.
Then
\begin{equation}
\Re\bigg(\bigg\langle
\ln\bigg(\frac{G_{\gg+1}(z(P),\cdot)-y(P)}{G_{\gg+1}(z(P),\cdot)+y(P)}\bigg)\bigg\rangle\bigg)=
0, \quad P\in\Gamma, \lb{4.22}
\end{equation}
where the closed simple curve $\Gamma\subset\calK_\gg$ denotes an
appropriate lift of $\gamma$ to $\calK_\gg$, yields the
contradiction
\begin{equation}
\Re\bigg(\bigg\langle
\ln\bigg(\frac{G_{\gg+1}(z(P),\cdot)-y(P)}{G_{\gg+1}(z(P),\cdot)+y(P)}\bigg)\bigg\rangle\bigg)=
0 \, \text{ for all $P$ in the interior of $\Gamma$} \lb{4.23}
\end{equation}
by Corollary 8.2.5 in \cite{Be86}. Therefore, since there are no
closed loops in $\sigma(H)$ and no analytic arc tends to infinity,
the resolvent set of $H$ is connected and hence path-connected,
proving $(v)$.
\end{proof}

\begin{remark} \lb{r4.2}
Here $\sigma\subset\bbC$ is called an {\it arc} if there exists a
parameterization $\gamma\in C([0,1])$ such that
$\sigma=\{\gamma(t)\,|\, t\in [0,1]\}$. The arc $\sigma$ is called
{\it simple} if there exists a parameterization $\gamma$ such that
$\gamma\colon [0,1]\to\bbC$ is injective.
\end{remark}

\appendix
\section{Hyperelliptic Curves and their Theta Functions} \lb{sA}
\renewcommand{\theequation}{A.\arabic{equation}}
\renewcommand{\thetheorem}{A.\arabic{theorem}}
\setcounter{theorem}{0} \setcounter{equation}{0}

We provide a brief summary of some of the fundamental notations
needed {}from the theory of hyperelliptic Riemann surfaces. More
details can be found in some of the standard textbooks \cite{FK92}
and \cite{Mu84} as well as in monographs dedicated to integrable
systems such as \cite[Ch.\ 2]{BBEIM94}, \cite[App.\ A, B]{GH03}.

Fix $\gg \in \bbN$. We intend to describe the hyperelliptic
Riemann surface $\calK_\gg$ of genus $\gg$ of the Toda-type curve
\eqref{1.2.19}, associated with the polynomial
\begin{align}
\begin{split}
&\calF_\gg(z,y)=y^2-R_{2\gg+2}(z)=0, \lb{a1} \\
&R_{2\gg+2}(z)=\prod_{m=0}^{2\gg+1}(z-E_m), \quad
\{E_m\}_{m=0}^{2\gg+1}\subset\bbC.
\end{split}
\end{align}
To simplify the discussion we will assume that the affine part of
$\calK_\gg$ is nonsingular, that is, we assume that
\begin{equation}
E_m \neq E_{m'} \text{ for } m\neq m', \; m,m'=0,\dots,2\gg+1
\lb{a2}
\end{equation}
throughout this appendix. Next we introduce an appropriate set of
(nonintersecting) cuts $\calC_j$ joining $E_{m(j)}$ and
$E_{m^\prime(j)}$, $j=1,\dots,\gg+1$, and denote
\begin{equation}
\calC=\bigcup_{j=1}^{\gg+1} \calC_j, \quad
\calC_j\cap\calC_k=\emptyset, \quad j\neq k.\lb{a3}
\end{equation}
Define the cut plane
\begin{equation}
\Pi=\bbC\backslash\calC, \lb{a4}
\end{equation}
and introduce the holomorphic function
\begin{equation}
R_{2\gg+2}(\cdot)^{1/2}\colon \Pi\to\bbC, \quad z\mapsto
\bigg(\prod_{m=0}^{2\gg+1}(z-E_m) \bigg)^{1/2}\lb{a5}
\end{equation}
on $\Pi$ with an appropriate choice of the square root branch in
\eqref{a5}. Next we define the set
\begin{equation}
\calM_{\gg}=\{(z,\sigma R_{2\gg+2}(z)^{1/2}) \mid z\in\bbC,\;
\sigma\in\{1,-1\} \}\cup \{P_{\infty_+},P_{\infty_-}\} \label{a6}
\end{equation}
by extending $R_{2\gg+2}(\cdot)^{1/2}$ to $\calC$. The
hyperelliptic curve $\calK_\gg$ is then the set $\calM_{\gg}$ with
its natural complex structure obtained upon gluing the two sheets
of $\calM_{\gg}$ crosswise along the cuts. Moreover, we introduce
the set of branch points
\begin{equation}
\calB(\calK_\gg)=\{(E_m,0)\}_{m=0}^{2\gg+1}. \lb{a7}
\end{equation}
Points $P\in\calK_\gg\backslash\{P_{\infty_{\pm}}\}$ are denoted
by
\begin{equation}
P=(z,\sigma R_{2\gg+2}(z)^{1/2})=(z,y),  \lb{a8}
\end{equation}
where $y(P)$ denotes the meromorphic function on $\calK_\gg$
satisfying $\calF_\gg(z,y)=y^2-R_{2\gg+2}(z)=0$ and
\begin{equation}
y(P)\underset{\zeta\to
0}{=}\mp\bigg(1-\f12\bigg(\sum_{m=0}^{2\gg+1}E_m\bigg)\zeta
+\Oh(\zeta^2)\bigg)\zeta^{-\gg-1} \text{  as $P\to
P_{\infty_\pm}$,} \; \zeta=1/z.     \lb{a55g}
\end{equation}

In addition, we introduce the holomorphic sheet exchange map
(involution)
\begin{equation}
*\colon\calK_\gg\to\calK_\gg, \quad P=(z,y)\mapsto P^*=(z,-y), \;
P_{\infty_\pm}\mapsto P^*_{\infty_\pm}=P_{\infty_\mp}  \lb{a9}
\end{equation}
and the two meromorphic projection maps
\begin{equation}
\tilde\pi\colon\calK_\gg\to\bbC\cup\{\infty\}, \quad
P=(z,y)\mapsto z, \; P_{\infty_\pm}\mapsto \infty \lb{a10}
\end{equation}
and
\begin{equation}
y\colon\calK_\gg\to\bbC\cup\{\infty\}, \quad P=(z,y)\mapsto y, \;
P_{\infty_\pm}\mapsto \infty.  \lb{a11}
\end{equation}
Thus the map $\tilde\pi$ has a pole of order 1 at $P_{\infty_\pm}$
and $y$ has a pole of order $\gg+1$ at $P_{\infty_\pm}$. Moreover,
\begin{equation}
\tilde\pi(P^*)=\tilde\pi(P), \quad y(P^*)=-y(P), \quad
P\in\calK_\gg. \lb{a12}
\end{equation}
As a result, $\calK_\gg$ is a two-sheeted branched covering of the
Riemann sphere $\bbC\bbP^1$ ($\cong\bbC\cup\{\infty\}$) branched
at the $2\gg+4$ points $\{(E_m,0)\}_{m=0}^{2\gg+1},
P_{\infty_\pm}$. $\calK_\gg$ is compact since $\tilde\pi$ is open
and $\bbC\bbP^1$ is compact. Therefore, the compact hyperelliptic
Riemann surface resulting in this manner has topological genus
$\gg$.

Next we introduce the upper and lower sheets $\Pi_{\pm}$ by
\begin{equation}
\Pi_{\pm}=\{(z,\pm  R_{2\gg+2}(z)^{1/2})\in \calM_\gg \mid
z\in\Pi\} \lb{a13}
\end{equation}
and the associated charts
\begin{equation}
\zeta_\pm \colon \Pi_\pm\to \Pi, \quad P\mapsto z.\lb{a14}
\end{equation}

Let $\{a_j,b_j\}_{j=1}^\gg$ be a homology basis for $\calK_\gg$
with intersection matrix of the cycles satisfying
\begin{equation}
a_j\circ b_k=\delta_{j,k}, \quad a_j\circ a_k=0, \quad b_j\circ
b_k=0, \quad j,k=1,\dots,\gg. \lb{a15}
\end{equation}
Associated with the homology basis $\{a_j,b_j\}_{j=1}^\gg$ we also
recall the canonical dissection of $\calK_\gg$ along its cycles
yielding the simply connected interior $\hatt \calK_ \gg$ of the
fundamental polygon $\partial {\hatt \calK}_\gg$ given by
\begin{equation}
\partial  {\hatt \calK}_\gg =a_1 b_1 a_1^{-1} b_1^{-1}
a_2 b_2 a_2^{-1} b_2^{-1} \cdots a_\gg^{-1} b_\gg^{-1}. \lb{a16}
\end{equation}

Let $\calM (\calK_\gg)$ and $\calM^1 (\calK_\gg)$ denote the set
of meromorphic functions (0-forms) and meromorphic differentials
(1-forms) on $\calK_\gg$, respectively. The residue of a
meromorphic differential $\nu\in \calM^1 (\calK_\gg)$ at a point
$Q \in \calK_\gg$ is defined by
\begin{equation}
\text{res}_{Q}(\nu) =\frac{1}{2\pi i} \int_{\gamma_{Q}} \nu,
\lb{a17}
\end{equation}
where $\gamma_{Q}$ is a counterclockwise oriented smooth simple
closed contour encircling $Q$ but no other pole of $\nu$.
Holomorphic differentials are also called Abelian differentials of
the first kind. Abelian differentials of the second kind
$\omega^{(2)} \in \calM^1 (\calK_\gg)$ are characterized by the
property that all their residues vanish. Any meromorphic differential
$\omega^{(3)}$ on $\calK_\gg$ not of the first or second kind is said to be
of the third kind. A differential of the third kind $\omega^{(3)} \in
\calM^1 (\calK_\gg)$ is usually normalized by vanishing of its
$a$-periods, that is,
\begin{equation}
\int_{a_j} \omega^{(3)} =0, \quad  j=1,\dots,\gg. \lb{a20}
\end{equation}
A normal differential $\omega_{P_1,P_2}^{(3)}$, associated with
two distinct points $P_1, \, P_2 \in \hat{\calK}_\gg$, by
definition, has simple poles at $P_1$ and $P_2$ with residues $+1$
at $P_1$ and $-1$ at $P_2$ and vanishing $a$-periods. If
$\omega_{P,Q}^{(3)}$ is a normal differential of the third kind
associated with $P, \, Q \in \hat{\calK}_\gg$, holomorphic on
$\calK_\gg \backslash \{P,Q\}$, then
\begin{equation}
\int_{b_j} \omega_{P,Q}^{(3)} =2\pi i \int_P^Q \omega_j, \quad
j=1,\dots,\gg. \lb{a21}
\end{equation}
We shall always assume (without loss of generality) that all poles
of $\omega^{(3)}$ on $\calK_\gg$ lie on
$\hat{\calK}_\gg$ (i.e., not on $\partial\hat{\calK}_\gg$).

Using local charts one infers that $d z/y$ is a holomorphic
differential on $\calK_\gg$ with zeros of order $\gg-1$ at
$P_{\infty_\pm}$ and hence
\begin{equation}
\eta_j=\frac{z^{j-1}d z}{y}, \quad j=1,\dots,\gg, \lb{a22}
\end{equation}
form a basis for the space of holomorphic differentials on
$\calK_\gg$. Introducing the invertible matrix $C$ in $\bbC^\gg$
\begin{align}
C & =\big(C_{j,k}\big)_{j,k=1,\dots,\gg}, \quad C_{j,k}
= \int_{a_k} \eta_j, \lb{a23} \\
\underline{c} (k) & = (c_1(k), \dots, c_\gg(k)), \quad c_j (k)
=\big(C^{-1}\big)_{j,k}, \quad j,k=1,\dots,\gg, \lb{a24}
\end{align}
the normalized differentials $\ome_j$ for $j=1,\dots,\gg$,
\begin{equation}
\ome_j = \sum_{\ell=1}^\gg c_j (\ell) \eta_\ell, \quad \int_{a_k}
\ome_j = \delta_{j,k}, \quad j,k=1,\dots,\gg, \lb{a25}
\end{equation}
form a canonical basis for the space of holomorphic differentials
on $\calK_\gg$.

In the chart $(U_{P_{\infty_\pm}}, \zeta_{P_{\infty_\pm}})$
induced by $1/\tilde\pi$ near $P_{\infty_\pm}$ one infers,
\begin{align}
{\ul \omega} & = (\omega_1,\dots,\omega_\gg)=
     \mp \sum_{j=1}^\gg \f{\uc (j)
\zeta^{\gg-j}d\zeta}{\big(\prod_{m=0}^{2\gg+1}
(1-\zeta E_m) \big)^{1/2}} \lb{a26} \\
& = \pm \bigg( \uc (\gg) +\zeta\big[ \frac12 \uc (\gg)
\sum_{m=0}^{2\gg+1} E_m +\uc (\gg-1) \big]  + \Oh(\zeta^2) \bigg)
d\zeta \text{ as $P\to P_{\infty_\pm}$,}
\no \\
& \hspace*{9.155cm} \zeta=1/z. \no
\end{align}

The matrix $\tau=\big(\tau_{j,\ell}\big)_{j,\ell=1}^\gg$  in
$\bbC^{\gg\times\gg}$ of $b$-periods defined by
\begin{equation}
\tau_{j,\ell}=\int_{b_j}\omega_\ell, \quad j,\ell=1, \dots,\gg
\label{a28}
\end{equation}
satisfies
\begin{equation}
\Im(\tau)>0 \, \text{ and } \, \tau_{j,\ell}=\tau_{\ell,j},
\;\, j,\ell =1,\dots,\gg.  \lb{a29}
\end{equation}

Associated with the matrix $\tau$ one introduces the period
lattice
\begin{equation}
L_\gg = \{ \ul z \in\bbC^\gg \mid \ul z = \ul m + \ul n\tau, \;
\ul m, \ul n \in\bbZ^\gg\} \lb{a30}
\end{equation}
and the Riemann theta function associated with $\calK_\gg$ and the
given homology basis $\{a_j,b_j\}_{j=1,\dots,\gg}$,
\begin{equation}
\theta(\ul z)=\sum_{\ul n\in\bbZ^\gg}\exp\big(2\pi i(\ul n,\ul
z)+\pi i(\ul n, \ul n\tau)\big), \quad \ul z\in\bbC^\gg,
\label{a31}
\end{equation}
where $(\ul u, \ul v)= \ol{\ul u}\,\ul v^\top =\sum_{j=1}^\gg
\overline{u}_j v_j$ denotes the scalar product in $\bbC^\gg$. It
has the following fundamental properties
\begin{align}
& \theta(z_1, \ldots, z_{j-1}, -z_j, z_{j+1}, \ldots, z_n) =\theta
(\ul z), \lb{a32}\\
& \theta (\ul z +\ul m + \ul n\tau) =\exp \big(-2 \pi i (\ul n,\ul
z) -\pi i (\ul n, \ul n\tau) \big) \theta (\ul z), \quad \ul m,
\ul n \in\bbZ^\gg. \lb{a33}
\end{align}

Next we briefly describe some consequences of a change of homology
basis. Let
\begin{equation}
\{a_1,\dots,a_\gg,b_1,\dots,b_\gg\} \lb{a34}
\end{equation}
be a canonical homology basis on $\calK_\gg$ with intersection
matrix satisfying \eqref{a15} and let
\begin{equation}
\{a'_1,\dots,a'_\gg,b'_1,\dots,b'_\gg\} \lb{a35}
\end{equation}
be a homology basis on $\calK_\gg$ related to \eqref{a34} by
\begin{equation}
\begin{pmatrix} {\ul a'}^\top \\ {\ul b'}^\top \end{pmatrix}
= X \begin{pmatrix} \ul a^\top \\ \ul b^\top \end{pmatrix},
\lb{a36}
\end{equation}
where
\begin{align}
\ul a^\top &=(a_1,\dots,a_\gg)^\top, \;\;\;\;\, \ul b^\top
=(b_1,\dots,b_\gg)^\top, \no \\
{\ul a'}^\top &=(a'_1,\dots,a'_\gg)^\top, \quad
{\ul b'}^\top =(b'_1,\dots,b'_\gg)^\top, \lb{a37} \\
X&=\begin{pmatrix} A & B \\ C & D \end{pmatrix} \lb{a38}
\end{align}
with $A,B,C$, and $D$ being $\gg\times \gg$ matrices with integer
entries. Then \eqref{a35} is also a canonical homology basis on
$\calK_\gg$ with intersection matrix satisfying \eqref{a15} if and
only if
\begin{equation}
X \in \Sp(\gg,\bbZ), \lb{a39}
\end{equation}
where
\begin{equation}
\Sp(\gg,\bbZ)=\left\{X=\begin{pmatrix} A & B \\ C & D
\end{pmatrix}\,\bigg|\, X\begin{pmatrix} 0 & I_\gg
\\  -I_\gg & 0 \end{pmatrix}X^\top=\begin{pmatrix} 0 &
I_\gg \\ -I_\gg & 0 \end{pmatrix}, \, \det(X)=1\right\} \lb{a40}
\end{equation}
denotes the symplectic modular group (here  $A,B,C,D$ in $X$ are
again $\gg\times \gg$ matrices with integer entries). If
$\{\omega_j\}_{j=1}^\gg$ and $\{\omega'_j\}_{j=1}^\gg$ are the
normalized bases of holomorphic differentials corresponding to the
canonical homology bases \eqref{a34} and \eqref{a35}, with $\tau$
and $\tau'$ the associated $b$ and $b'$-periods of
${\ul\omega}=\omega_1,\dots,\omega_\gg$ and
${\ul\omega'}=\omega'_1,\dots,\omega'_\gg$, respectively, then one
computes
\begin{equation}
{\ul \omega'}={\ul\omega}(A+B\tau)^{-1}, \quad
\tau'=(C+D\tau)(A+B\tau)^{-1}. \lb{a41}
\end{equation}

Fixing a base point $Q_0\in\calK_\gg\backslash\{P_{\infty_\pm}\}$,
one denotes by $J(\calK_\gg) = \bbC^\gg/L_\gg$ the Jacobi variety
of $\calK_\gg$, and defines the Abel map $\underline{A}_{Q_0}$ by
\begin{equation}
\underline{A}_{Q_0} \colon \calK_n \to J(\calK_\gg), \quad
\underline{A}_{Q_0}(P)= \bigg(\int_{Q_0}^P
\omega_1,\dots,\int_{Q_0}^P \omega_\gg \bigg) \pmod{L_\gg}, \quad
P\in\calK_\gg. \label{a42}
\end{equation}

Next, consider the vector $\ul{U}_{0}^{(3)}$ of $b$-periods of
$\omega_{P_{\infty_+},{P_{\infty_-}}}^{(3)}/(2\pi i)$, the
normalized differential of the third kind, holomorphic on
$\calK_\gg\backslash\{P_{\infty_\pm}\}$,
\begin{equation}
\ul{U}_{0}^{(3)}=\big({U}_{0,1}^{(3)},\dots,{U}_{0,\gg}^{(3)}\big),
\quad {U}_{0,j}^{(3)}=\f{1}{2\pi i}\int_{b_j}
\omega_{P_{\infty_+},{P_{\infty_-}}}^{(3)}, \;\, j=1,\dots,\gg.
\lb{a27}
\end{equation}
Using \eqref{a21} one then computes 
\begin{equation}
{\ul U}_{0}^{(3)}=\ul{A}_{P_{\infty -}}({P_{\infty
+}})=2\ul{A}_{Q_{0}}({P_{\infty +}}), \lb{a27a}
\end{equation}
where $Q_0$ is chosen to be a branch point of $\calK_\gg$, $Q_0\in
\calB(\calK_\gg)$, in the last part of \eqref{a27a}.

Similarly, one introduces
\begin{equation}
\ul \alpha_{Q_0}  \colon \Div(\calK_\gg) \to J(\calK_\gg),\quad
\calD \mapsto \ul \alpha_{Q_0} (\calD) =\sum_{P \in \calK_\gg}
\calD (P) \ul A_{Q_0} (P), \label{a43}
\end{equation}
where $\Div(\calK_\gg)$ denotes the set of divisors on
$\calK_\gg$. Here a map $\calD \colon \calK_\gg \to \bbZ$ is
called a divisor on $\calK_\gg$ if $\calD(P)\neq0$ for only
finitely many $P\in\calK_\gg$. (In the main body of this paper we
will choose $Q_0$ to be one of the branch points, i.e.,
$Q_0\in\calB(\calK_\gg)$, and we will always choose
the same path of integration from $Q_0$ to $P$ in all Abelian
integrals.) For subsequent use in Remark \ref{raa26a} we also
introduce
\begin{align}
\hua_{Q_0} & \colon\hatt{\calK}_\gg\to\bbC^\gg, \lb{a44} \\
& \quad \, P\mapsto\hua_{Q_0}(P) =\big(\hatt A_{Q_0,1}(P),\dots,\hatt
A_{Q_0,\gg}(P)\big)
=\bigg(\int_{Q_0}^P\omega_1,\dots,\int_{Q_0}^P\omega_\gg\bigg) \no
\end{align}
and
\begin{equation}
\hatt {\ul \al}_{Q_0}  \colon \Div(\hatt\calK_\gg) \to \bbC^\gg,
\quad \calD \mapsto \hatt {\ul \al}_{Q_0} (\calD) =\sum_{P \in
\hatt\calK_\gg} \calD (P) \hua_{Q_0} (P). \lb{a45}
\end{equation}

In connection with divisors on $\calK_\gg$ we will employ the
following (additive) notation,
\begin{align}
&\calD_{Q_0\ul Q}=\calD_{Q_0}+\calD_{\ul Q}, \quad \calD_{\ul
Q}=\calD_{Q_1}+\cdots +\calD_{Q_m}, \lb{a46} \\
& {\ul Q}=\{Q_1, \dots ,Q_m\} \in \sym^m \calK_\gg, \quad
Q_0\in\calK_\gg, \;\, m\in\bbN, \no
\end{align}
where for any $Q\in\calK_\gg$,
\begin{equation} \lb{a47}
\calD_Q \colon  \calK_\gg \to\bbN_0, \quad P \mapsto  \calD_Q (P)=
\begin{cases} 1 & \text{for $P=Q$},\\
0 & \text{for $P\in \calK_\gg\backslash \{Q\}$}, \end{cases}
\end{equation}
and $\sym^m \calK_\gg$ denotes the $m$th symmetric product of
$\calK_\gg$. In particular, $\sym^m \calK_\gg$ can be identified
with the set of nonnegative divisors $0 \leq \calD \in
\Div(\calK_\gg)$ of degree $m\in\bbN$.

For $f\in \calM (\calK_\gg) \backslash \{0\}$, $\omega \in \calM^1
(\calK_\gg) \backslash \{0\}$ the divisors of $f$ and $\omega$ are
denoted by $(f)$ and $(\omega)$, respectively.  Two divisors
$\calD$, $\calE\in \Div(\calK_\gg)$ are called equivalent, denoted
by $\calD \sim \calE$, if and only if $\calD -\calE =(f)$ for some
$f\in\calM (\calK_\gg) \backslash \{0\}$.  The divisor class
$[\calD]$ of $\calD$ is then given by $[\calD] =\{\calE \in
\Div(\calK_{\gg})\mid\calE \sim \calD\}$.  We recall that
\begin{equation}
\deg ((f))=0, \quad \deg ((\omega)) =2(\gg-1), \quad f\in\calM (\calK_\gg)
\backslash \{0\}, \;  \omega\in \calM^1 (\calK_\gg) \backslash
\{0\}, \lb{a48}
\end{equation}
where the degree $\deg (\calD)$ of $\calD$ is given by $\deg
(\calD) =\sum_{P\in \calK_\gg} \calD (P)$.  It is customary to
call $(f)$ (respectively, $(\omega)$) a principal (respectively,
canonical) divisor.

Introducing the complex linear spaces
\begin{align}
\calL (\calD) & =\{f\in \calM (\calK_\gg)\mid f=0
      \text{ or } (f) \geq \calD\}, \quad 
r(\calD) =\dim_\bbC \calL (\calD),
\lb{a49}\\
\calL^1 (\calD) & =
      \{ \omega\in \calM^1 (\calK_\gg)\mid \omega=0
      \text{ or } (\omega) \geq
\calD\},\quad  i(\calD) =\dim_\bbC \calL^1 (\calD)  \lb{a50}
\end{align}
(with $i(\calD)$ the index of specialty of $\calD$), one infers
that $\deg (\calD)$, $r(\calD)$, and $i(\calD)$ only depend on the
divisor class $[\calD]$ of $\calD$.  Moreover, we recall the
following fundamental facts.

\begin{theorem} \lb{thm1}
Let $\calD \in \Div(\calK_\gg)$, $\omega \in \calM^1 (\calK_\gg)
\backslash \{0\}$. Then,
\begin{equation}
      i(\calD) =r(\calD-(\omega)), \quad \gg\in\bbN_0.
\lb{a51}
\end{equation}
The Riemann-Roch theorem reads
\begin{equation}
r(-\calD) =\deg (\calD) + i (\calD) -\gg+1, \quad n\in\bbN_0.
\lb{a52}
\end{equation}
By Abel's theorem, $\calD\in \Div(\calK_\gg)$, $\gg\in\bbN$, is
principal if and only if
\begin{equation}
\deg (\calD) =0 \text{ and } \ul \alpha_{Q_0} (\calD) =\ul{0}.
\lb{a53}
\end{equation}
Finally, assume $\gg\in\bbN$. Then $\ul \alpha_{Q_0} :
\Div(\calK_\gg) \to J(\calK_\gg)$ is surjective $($Jacobi's
inversion theorem$)$.
\end{theorem}

\begin{theorem} \lb{thm2}
Let $\calD_{\ul Q} \in \sym^{\gg} \calK_\gg$, $\ul Q=\{Q_1,
\ldots, Q_\gg\}$.  Then,
\begin{equation}
1 \leq i (\calD_{\ul Q} ) =s \lb{a54}
\end{equation}
if and only if there are $s$ pairs of the type $\{P,
P^*\}\subseteq \{Q_1,\ldots, Q_\gg\}$ $($this includes, of course,
branch points for which $P=P^*$$)$. One has $s\leq
\gg/2$.
\end{theorem}

Next, we denote by $\ul \Xi_{Q_0}=(\Xi_{Q_{0,1}}, \dots,
\Xi_{Q_{0,\gg}})$ the vector of Riemann constants,
\begin{equation}
\Xi_{Q_{0,j}}=\frac12(1+\tau_{j,j})- \sum_{\substack{\ell=1 \\
\ell\neq j}}^\gg\int_{a_\ell} \omega_\ell(P)\int_{Q_0}^P\omega_j,
\quad j=1,\dots,\gg. \lb{a55}
\end{equation}

\begin{theorem} \lb{thm3}
Let $\ul Q =\{Q_1,\dots,Q_\gg\}\in \sym^{\gg} \calK_\gg$ and
assume $\calD_{\ul Q}$ to be nonspecial, that is, $i(\calD_{\ul
Q})=0$. Then,
\begin{equation}
\theta\big(\ul {\Xi}_{Q_0} -\ul {A}_{Q_0}(P) + \alpha_{Q_0}
(\calD_{\ul Q})\big)=0 \text{ if and only if }
P\in\{Q_1,\dots,Q_\gg\}. \lb{a56}
\end{equation}
\end{theorem}

\begin{remark} \lb{raa26a}
In Section \ref{s2} we dealt with theta function expressions of
the type
\begin{equation}
\psi(P)=\f{\theta(\ul\Xi_{Q_0}-\ul
A_{Q_0}(P)+\ul\alpha_{Q_0}(\calD_1))} {\theta(\ul\Xi_{Q_0}-\ul
A_{Q_0}(P)+\ul\alpha_{Q_0}(\calD_2))} \exp\bigg(-c \int_{Q_0}^P
\Omega^{(3)}\bigg), \quad P\in\calK_\gg, \lb{a57}
\end{equation}
where $\calD_j\in\sym^\gg\calK_\gg$, $j=1,2$, are nonspecial
positive divisors of degree $\gg$, $c\in\bbC$ is a constant, and
$\Omega^{(3)}$ is a normalized differential of the third kind with
a prescribed singularity at $P_{\infty_\pm}$. Even though we agree
to always choose identical paths of integration {}from $P_0$ to
$P$ in all Abelian integrals \eqref{a57}, this is not sufficient
to render $\psi$ single-valued on $\calK_\gg$. To achieve
single-valuedness one needs to replace $\calK_\gg$ by its simply
connected canonical dissection $\hatt\calK_{\gg}$ and then replace
$\ul A_{Q_0}$ and $\ul \alpha_{Q_0}$ in \eqref{a57} with
${\hua}_{Q_0}$ and $\hatt{\ul \alpha}_{Q_0}$ as introduced in
\eqref{a44} and \eqref{a45}. In particular, one regards $a_j,b_j$,
$j=1,\dots,\gg$, as curves (being a part of
$\partial\hatt\calK_\gg$, cf. \eqref{a16}) and not as homology
classes. Similarly, one then replaces $\uxi_{Q_0}$ by \,$\hatt
\uxi_{Q_0}$ (replacing $\ul A_{Q_0}$ by ${\hua}_{Q_0}$ in
\eqref{a55}, etc.). Moreover, in connection with $\psi$, one
introduces the vector of $b$-periods $\ul U^{(3)}$ of
$\Omega^{(3)}$ by
\begin{equation}
\ul U^{(3)}=(U_1^{(3)},\dots,U_{\gg}^{(3)}), \quad
U_j^{(3)}=\f{1}{2\pi i}\int_{b_j} \Omega^{(3)}, \quad
j=1,\dots,\gg, \lb{a58}
\end{equation}
and then renders $\psi$ single-valued on $\hatt\calK_\gg$ by
requiring
\begin{equation}
\hatt{\ul\alpha}_{Q_0}(\calD_1)-\hatt{\ul\alpha}_{Q_0}(\calD_2) =c
\,\ul U^{(3)} \lb{a59}
\end{equation}
$($as opposed to merely
$\ul\alpha_{Q_0}(\calD_1)-\ul\alpha_{Q_0}(\calD_2)=c \,\ul U^{(3)}
\pmod {L_\gg}$$)$. Actually, by \eqref{a33},
\begin{equation}
\hatt{\ul\alpha}_{Q_0}(\calD_1)-\hatt{\ul\alpha}_{Q_0}(\calD_2) -
c\, \ul U^{(3)}\in\bbZ^\gg, \lb{a60}
\end{equation}
suffices to guarantee single-valuedness of $\psi$ on
$\hatt\calK_\gg$. Without the replacement of $\ul A_{Q_0}$ and
$\ul \alpha_{Q_0}$ by ${\hua}_{Q_0}$ and $\hatt{\ul \alpha}_{Q_0}$
in \eqref{a57} and without the assumption \eqref{a59} $($or
\eqref{a60}$)$, $\psi$ is a multiplicative $($multi-valued$)$
function on $\calK_\gg$, and then most effectively discussed by
introducing the notion of characters on $\calK_\gg$ $($cf.\
\cite[Sect.\ III.9]{FK92}$)$. For simplicity, we decided to avoid
the latter possibility and throughout this paper will always
tacitly assume \eqref{a59} or \eqref{a60}.
\end{remark}

\section{Restrictions on $\ul B = \ul U_0^{(3)}$} \lb{sB}
\renewcommand{\theequation}{B.\arabic{equation}}
\renewcommand{\thetheorem}{B.\arabic{theorem}}
\setcounter{theorem}{0} \setcounter{equation}{0}

The purpose of this appendix is to prove the result \eqref{2.66},
$\ul B=\ul U^{(3)}_0 \in \bbR^\gg$, for some choice of homology
basis $\{a_j, b_j\}_{j=1}^\gg$ on $\calK_\gg$ as recorded in
Remark \ref{r2.8}.

To this end we first recall a few notions in connection with
periodic meromorphic functions of $p$ complex variables.

\begin{definition} \lb{dB.1}
Let $p\in\bbN$ and $F\colon\bbC^p\to\bbC\cup\{\infty\}$ be
meromorphic (i.e., a ratio of two entire functions of $p$ complex
variables).
Then, \\
(i) $\ul \omega=(\omega_1,\dots,\omega_p)\in\bbC^p\backslash\{0\}$
is called a {\it period} of $F$ if
\begin{equation}
F(\ul z+\ul \omega)=F(\ul z) \lb{B.1}
\end{equation}
for all $\ul z\in\bbC^p$ for which $F$ is analytic. The set of all
periods of $F$ is denoted by $\calP_F$. \\
(ii) $F$ is called {\it degenerate} if it depends on less than $p$
complex variables; otherwise, $F$ is called {\it nondegenerate}.
\end{definition}

\begin{theorem} \lb{tB.2}
Let $p\in\bbN$, $F\colon\bbC^p\to\bbC\cup\{\infty\}$ be
meromorphic, and
$\calP_F$ be the set of all periods of $F$. Then either \\
$(i)$ $\calP_F$ has a finite limit point, \\
or \\
$(ii)$ $\calP_F$ has no finite limit point. \\
In case $(i)$, $\calP_F$ contains {\it infinitesimal periods}
$($i.e., sequences of nonzero periods converging to zero$)$. In
addition, in case $(i)$ each period is a limit point of periods
and hence $\calP_F$ is a perfect set. \\
Moreover, $F$ is degenerate if and only if $F$ admits
infinitesimal periods. In particular, for nondegenerate functions
$F$ only alternative $(ii)$ applies.
\end{theorem}

Next, let $\ul\omega_q\in\bbC^p\backslash\{0\}$, $q=1,\dots,r$ for
some $r\in\bbN$. Then $\ul\omega_1,\dots,\ul\omega_r$ are called
{\it linearly independent over $\bbZ$ $($resp.\ $\bbR$$)$} if
\begin{align}
&\nu_1\ul\omega_1+\cdots+\nu_r\ul\omega_r=0, \quad \nu_q\in\bbZ
\text{ (resp., $\nu_q\in\bbR$)}, \;\, q=1,\dots,r, \no \\
&\text{implies } \nu_1=\cdots=\nu_r=0. \lb{B.2}
\end{align}
Clearly, the maximal number of vectors in $\bbC^p$ linearly
independent over $\bbR$ equals $2p$.

\begin{theorem} \lb{tB.3}
Let $p\in\bbN$. \\
$(i)$ If $F\colon\bbC^p\to\bbC\cup\{\infty\}$ is a nondegenerate
meromorphic function with periods
$\ul\omega_q\in\bbC^p\backslash\{0\}$, $q=1,\dots,r$, $r\in\bbN$,
linearly independent over $\bbZ$, then
$\ul\omega_1,\dots,\ul\omega_r$ are also linearly independent over
$\bbR$. In particular, $r\leq 2p$.
\\
$(ii)$ A nondegenerate entire function $F\colon\bbC^p\to\bbC$
cannot have more than $p$ periods linearly independent over $\bbZ$
$($or $\bbR$$)$.
\end{theorem}

For $p=1$, $\exp(z)$, $\sin(z)$ are examples of entire functions
with precisely one period. Any non-constant doubly periodic
meromorphic function of one complex variable is elliptic (and
hence indeed has poles).

\begin{definition} \lb{dB.4}
Let $p, r\in\bbN$. A system of periods
$\ul\omega_q\in\bbC^p\backslash\{0\}$, $q=1,\dots,r$, of a
nondegenerate meromorphic function
$F\colon\bbC^p\to\bbC\cup\{\infty\}$, linearly independent over
$\bbZ$, is called {\it fundamental} or a {\it basis} of periods
for $F$ if every period $\ul\omega$ of $F$ is of the form
\begin{equation}
\ul\omega =m_1\ul\omega_1+\cdots+m_r\ul\omega_r \, \text{ for some
$m_q\in\bbZ$, $q=1,\dots,r$.} \lb{B.3}
\end{equation}
\end{definition}

The representation of $\ul\omega$ in \eqref{B.3} is unique since
by hypothesis $\ul\omega_1,\dots,\ul\omega_r$ are linearly
independent over $\bbZ$. In addition, $\calP_F$ is countable in
this case. (This rules out case $(i)$ in Theorem \ref{tB.2} since
a perfect set is uncountable. Hence, one does not have to assume
that $F$ is nondegenerate in Definition \ref{dB.4}.)

This material is standard and can be found, for instance, in
\cite[Ch.\ 2]{Ma92}.

\vspace*{2mm}

Next, returning to the Riemann theta function $\theta(\ul\cdot)$
in \eqref{a31}, we introduce the vectors $\{\ul e_j\}_{j=1}^\gg,
\{\ul\tau_j\}_{j=1}^\gg \subset\bbC^\gg\backslash\{0\}$ by
\begin{equation}
\ul e_j = (0,\dots,0,\underbrace{1}_{j},0,\dots,0), \quad \ul
\tau_j = \ul e_j \tau, \quad j=1,\dots,\gg. \lb{B.4}
\end{equation}
Then
\begin{equation}
\{\ul e_j\}_{j=1}^\gg \lb{B.5}
\end{equation}
is a basis of periods for the entire (nondegenerate) function
$\theta(\ul\cdot)\colon\bbC^\gg\to\bbC$. Moreover, fixing
$k\in\{1,\dots,\gg\}$, then
\begin{equation}
\{\ul e_j, \ul\tau_j\}_{j=1}^\gg \lb{B.6}
\end{equation}
is a basis of periods for the meromorphic function
$\partial_{z_k}\ln\big(\f{\theta(\ul\cdot)}{\theta(\ul\cdot+\ul{V})}\big)\colon\bbC^\gg\to\bbC\cup\{\infty\}$,
$\ul{V}\in \bbC^{\gg}$ (cf.\ \eqref{a33} and \cite[p.\ 91]{FK92}).

Next, let $\ul A\in\bbC^\gg$, $\ul
D=(D_1,\dots,D_\gg)\in\bbR^\gg$, $D_j\in\bbR\backslash\{0\}$,
$j=1,\dots,\gg$, and consider
\begin{align}
\begin{split}
f_{k}\colon \bbR\to\bbC, \quad f_{k}(n)&=\partial_{z_k}
\ln\bigg(\f{\theta(\ul A-\ul z)}{\theta(\ul C -\ul
z)}\bigg)\bigg|_{\ul z=\ul D n} \lb{B.7} \\
&=\partial_{z_k} \ln\bigg(\f{\theta(\ul A-\ul z\diag(\ul
D))}{\theta(\ul C-\ul z\diag(\ul D))}\bigg)\bigg|_{\ul
z=(n,\dots,n)}.
\end{split}
\end{align}
Here $\diag(\ul D)$ denotes the diagonal matrix
\begin{equation}
\diag(\ul D)= \big(D_j\delta_{j,j'}\big)_{j,j'=1}^\gg. \lb{B.8}
\end{equation}
Then the quasi-periods $D_j^{-1}$, $j=1,\dots,\gg$, of $f_{k}$ are
in a one-to-one correspondence with the periods of
\begin{equation}
F_{k}\colon\bbC^\gg\to\bbC\cup\{\infty\}, \quad F_{k}(\ul
z)=\partial_{z_{k}} \ln\bigg(\f{\theta(\ul A-\ul z\diag(\ul
D))}{\theta(\ul C-\ul z\diag(\ul D))}\bigg) \lb{B.9}
\end{equation}
of the special type
\begin{equation}
\ul e_j \big(\diag (\ul D)\big)^{-1} =
\big(0,\dots,0,\underbrace{D_j^{-1}}_{j},0,\dots,0\big). \lb{B.10}
\end{equation}
Moreover,
\begin{equation}
f_{k}(n)=F_{k}(\ul z)|_{\ul z=(n,\dots,n)}, \quad n\in\bbZ.
\lb{B.11}
\end{equation}

\begin{theorem} \lb{tB.5}
Suppose $a$ and $b$  in \eqref{1.2.45} to be quasi-periodic. Then
there exists a homology basis $\{\ti a_j, \ti b_j\}_{j=1}^\gg$ on
$\calK_\gg$ such that the vector $\wti{\ul B}=\wti{\ul U}^{(3)}_0$
with $\wti{\ul U}^{(3)}_0$ the vector of $\ti b$-periods of the
corresponding normalized differential of the third kind, $\wti
\omega^{(3)}_{P_{\infty+},P_{\infty-}}$, satisfies the constraint
\begin{equation}
\wti {\ul B}=\wti{\ul U}^{(3)}_0 \in \bbR^\gg. \lb{B.12}
\end{equation}
\end{theorem}
\begin{proof}
By \eqref{a27}, the vector of $b$-periods $\ul U^{(3)}_0$
associated with a given homology basis $\{a_j, b_j\}_{j=1}^\gg$ on
$\calK_\gg$ and the normalized differential of the third kind,
$\omega^{(3)}_{P_{\infty+},P_{\infty-}}$, is continuous with
respect to $E_0,\dots,E_{2\gg+1}$. Hence, we may assume in the
following that
\begin{equation}
B_j\neq 0, \;\, j=1,\dots,\gg, \quad \ul B=(B_1,\dots,B_\gg)
\lb{B.13}
\end{equation}
by slightly altering $E_0,\dots,E_{2\gg+1}$, if necessary. Using
\eqref{1.3.IM}, we may write
\begin{align}
\begin{split}
b(n)& =\Lambda_0
-\sum_{j=1}^{\gg}c_j(\gg)\frac{\partial}{\partial\omega_j}\ln\bigg(\f{\theta(\underline{\omega}+\ul
A -\ul B n)}{\theta(\underline{\omega}+\ul C -\ul B
n)}\bigg)\bigg|_{\underline{\omega}=0}
\\
&=\Lambda_0 -\sum_{j=1}^\gg c_j(\gg)
\partial_{z_{j}}
\ln\bigg(\f{\theta(\ul A-\ul z)}{\theta(\ul C -\ul
z)}\bigg)\bigg|_{\ul z=\ul B n}, \lb{B.14}
\end{split}
\end{align}
where by \eqref{1.3.IMB},
\begin{equation}
\ul{B}=\ul{U}_{0}^{(3)}. \lb{B.14a}
\end{equation}
Introducing the meromorphic (nondegenerate) function
$\calV\colon\bbC^\gg\to\bbC\cup\{\infty\}$ by
\begin{equation}
\calV(\ul z)=\Lambda_0 -\sum_{j=1}^n c_j(\gg)
\partial_{z_{j}}
\ln\bigg(\f{\theta(\ul A-\ul z\diag(\ul B))}{\theta(\ul C-\ul
z\diag(\ul B))}\bigg), \lb{B.15}
\end{equation}
one observes that
\begin{equation}
b(n)=\calV(\ul z)|_{\ul z=(n,\dots,n)}. \lb{B.16}
\end{equation}
In addition, $\calV$ has a basis of periods
\begin{equation}
\Big\{\ul e_j \big(\diag(\ul B)\big)^{-1}, \ul\tau_j
\big(\diag(\ul B)\big)^{-1}\Big\}_{j=1}^\gg \lb{B.17}
\end{equation}
by \eqref{B.6}, where
\begin{align}
\ul e_j \big(\diag(\ul B)\big)^{-1} &=
\big(0,\dots,0,\underbrace{B_j^{-1}}_{j},0,\dots,0\big),
\quad j=1,\dots,\gg, \lb{B.18} \\
\ul\tau_j\big(\diag(\ul B)\big)^{-1}
&=\big(\tau_{j,1}B_1^{-1},\dots, \tau_{j,\gg}B_\gg^{-1}\big),
\quad j=1,\dots,\gg. \lb{B.19}
\end{align}

By hypothesis, $b$ in \eqref{B.14} is quasi-periodic and hence has
$\gg$ real (scalar) quasi-periods. The latter are not necessarily
linearly independent over $\bbQ$ from the outset, but by slightly
changing the locations of branchpoints $\{E_m\}_{m=0}^{2\gg+1}$
into, say, $\{\wti E_m\}_{m=0}^{2\gg+1}$, one can assume they are.
In particular, since the period vectors in \eqref{B.17} are
linearly independent and the (scalar) quasi-periods of $b$ are in
a one-one correspondence with vector periods of $\calV$ of the
special form \eqref{B.18} (cf.\ \eqref{B.9}, \eqref{B.10}), there
exists a homology basis $\{\ti a_j, \ti b_j\}_{j=1}^\gg$ on
$\calK_\gg$ such that the vector $\ul {\wti B}=\wti{\ul
U}^{(3)}_0$, corresponding to the normalized differential of the
third kind, $\wti \omega^{(3)}_{P_{\infty+},P_{\infty-}}$ and this
particular homology basis, is real-valued. By continuity of
$\wti{\ul U}^{3}_0$ with respect to $\wti E_0,\dots,\wti
E_{2\gg+1}$, this proves \eqref{B.12}.
\end{proof}

{\bf Acknowledgments.}
We are grateful to Leonid Golinskii for pointing out references \cite{Na62}
and \cite{Na64} to us.



\begin{thebibliography}{10}
%
\bi{Ap86} A.\ I.\ Aptekarev, {\it Asymptotic properties of polynomials
orthogonal on a system of contours, and periodic motions of Toda lattices},
Math. USSR Sbornik {\bf 53}, 233--260 (1986).
%
\bi{BG05} V.\ Batchenko and F.\ Gesztesy, {\it On the spectrum of
Schr\"odinger operators with quasi-periodic algebro-geometric KdV
potentials}, J. Analyse Math. {\bf 95}, 333--387 (2005). 
%
\bi{Be86} A.\ F.\ Beardon, {\it A Primer on Riemann Surfaces},
London Math. Soc. Lecture Notes {\bf 78}, Cambridge University
Press, Cambridge, 1986.
%
\bibitem{BBEIM94} E.\ D.\ Belokolos, A.\ I.\ Bobenko, V.\ Z.\ Enol'skii,
A.\ R.\ Its, and V.\ B.\ Matveev, {\it Algebro-Geometric Approach to
Nonlinear Integrable Equations}, Springer, Berlin, 1994.
%
\bi{Be54} A.\ S.\ Besicovitch, {\it Almost Periodic Functions},
Dover, New York, 1954.
%
\bi{Bo47} H.\ Bohr, {\it Almost Periodic Functions}, Chelsea, New
York,1947.
%
\bi{BGHT98} W.\ Bulla F.\ Gesztesy, H.\ Holden, and G.\ Teschl.
{\it Algebro-geometric quasi-periodic finite-gap solutions of the
Toda and Kac-van Moerbeke hierarchies}, Memoirs  Amer. Math. Soc. {\bf 135},
No.\ 641, 1--79 (1998).
%
\bi{CK87} R.\ Carmona and S.\ Kotani, {\it Inverse spectral theory for
random Jacobi matrices},  J. Statist. Phys. {\bf 46}, 1091--1114 (1987).
%
\bi{CL90} R.\ Carmona and J.\ Lacroix, {\it Spectral Theory of
Random Schr\"odinger Operators}, Birkh\"auser, Boston, 1990.
%
\bi{CE70} R.\ S.\ Chisholm and W.\ N.\ Everitt, {\it On bounded
integral operators in the space of integrable-square functions},
Proc. Roy. Soc. Edinburgh Sect. A {\bf 69}, 199--204 (1970/71).
%
\bi{Co89} C.\ Corduneanu, {\it Almost Periodic Functions}, 2nd
ed., Chelsea, New York, 1989.
%
\bi{DT176} E.\ Date and S.\ Tanaka, {\it Analogue of inverse
scattering theory for the discrete {H}ill's equation and exact
solutions for the periodic {T}oda lattice}, Progr. Theoret. Phys.
{\bf 56}, 457--465 (1976).
%
\bi{DT276} E.\ Date and S.\ Tanaka, {\it Periodic multi-soliton
solutions of {K}orteweg--de {V}ries equation and {T}oda lattice},
Progr. Theoret. Phys. Suppl. {\bf 59}, 107--125 (1976).
%
\bi{DS83} F.\ Delyon and B.\ Souillard, {\it The rotation number for
finite difference operators and its properties}, Commun. Math. Phys. {\bf
89}, 415--426 (1983).
%
\bi{Du81} B.\ A.\ Dubrovin, {\it Theta functions and non-linear
equations}, Russian Math. Surveys {\bf 36:2}, 11--92 (1981).
%
\bi{DKN90} B.\ A.\ Dubrovin and I.\ M.\ Krichever and S.\ P.\
Novikov, {\it Integrable systems. {I}}, Dynamical Systems IV,
Springer, Berlin, 1990, pp.\ 173--280.
%
\bi{DMN76} B.\ A.\ Dubrovin, V.\ B.\ Matveev, and S.\ P.\ Novikov,
{\it Non-linear equations of Korteweg-de Vries type, finite-zone
linear operators, and Abelian varieties}, Russian Math. Surv. {\bf
31:1}, 59--146 (1976).
%
\bi{EE89} D.\ E.\ Edmunds and W.\ D.\ Evans, {\it Spectral Theory and
Differential Operators}, Clarendon Press, Oxford, 1989.
%
\bi{Ea67} M.\ S.\ P.\ Eastham, {\it Gaps in the essential spectrum
associated with singular differential operators}, Quart. J. Math.
{\bf 18}, 155--168 (1967).
%
\bi{Ea73} M.\ S.\ P.\ Eastham, {\it The Spectral Theory of
Periodic Differential Equations}, Scottish Academic Press,
Edinburgh and London, 1973.
%
\bi{FK92} H.\ M.\ Farkas and I.\ Kra, {\it Riemann Surfaces}, 2nd
ed., Springer, New York, 1992.
%
\bi{Fi74} A.\ M. Fink, {\it Almost Periodic Differential
Equations}, Lecture Notes in Math. {\bf 377}, Springer, Berlin,
1974.
%
\bi{Fl175} H.\ Flaschka, {\it Discrete and periodic illustrations
of some aspects of the inverse method}, in {\it Dynamical
{S}ystems, {T}heory and {A}pplications}, Lecture Notes in Phys.  
{\bf 38}, Springer, Berlin, 1975, pp.\ 441--466.
%
\bi{GH03} F.\ Gesztesy and H.\ Holden, {\it Soliton Equations and
Their Algebro-Geometric Solutions. Vol. I: $(1+1)$-Dimensional
Continuous Models},  Cambridge Studies in Advanced Mathematics,
Vol.\ 79, Cambridge Univ. Press, Cambridge, 2003.
%
\bi{GH05} F.\ Gesztesy and H.\ Holden, {\it Soliton Equations and
Their Algebro-Geometric Solutions. Vol. II: $(1+1)$-Dimensional
Discrete Models},  Cambridge Studies in Advanced Mathematics,
Cambridge Univ. Press, in preparation.
%
\bi{GHT05} F.\ Gesztesy, H.\ Holden, and G.\ Teschl, {\it The
algebro-geometric Toda hierarchy initial value problem for complex-valued
initial data}, in preparation.
%
\bi{Gl65} I.\ M. Glazman, {\it Direct Methods of Qualitative
Spectral Analysis of Singular Differential Operators}, Moscow,
1963. Engl. Transl. by Israel Program for Scientific
Translations, 1965.
%
\bi{Go85} S.\ Goldberg, {\it Unbounded Linear Operators}, Dover,
New York, 1985.
%
\bi{GK03} I.\ Ya.\ Goldsheid and B.\ A.\ Khoruzhenko, {\it Thouless
formula for random non-Hermitian Jacobi matrices}, preprint, 2003.
%
\bi{JM82} R.\ Johnson and J.\ Moser, {\it The rotation number for
almost periodic potentials}, Commun. Math. Phys. {\bf 84},
403--438 (1982).
%
\bi{Ka80} T.\ Kato, {\it Perturbation Theory for Linear
Operators}, corr. printing of the 2nd ed., Springer, Berlin, 1980.
%
\bi{KK03} E.\ Korotyaev and I.\ V.\ Krasovsky, {\it Spectral estimates
for periodic Jacobi matrices}, Comm. Math. Phys. {\bf 234}, 517--532
(2003).
%
\bi{Kr78} I.\ M.\ Krichever, {\it Algebraic curves and non-linear
difference equations}, Russian Math. Surveys {\bf 33:4}, 255--256
(1978).
%
\bi{Kr82} I.\ M.\ Krichever, {\it Algebro-geometric spectral
theory of the {S}chr\"odinger difference operator and the
{P}eierls model}, Soviet Math. Dokl. {\bf 26}, 194--198 (1982).
%
\bi{Kr182}I.\ M.\ Krichever, {\it The {P}eierls model}, Functional
Anal. Appl. {\bf 16}, 248--263 (1982).
%
\bi{Kr83} I.\ M.\ Krichever, {\it Nonlinear equations and elliptic
curves}, Revs. Sci. Tech., {\bf 23}, 51--90 (1983).
%
\bi{Le87} B.\ M.\ Levitan, {\it Inverse Sturm--Liouville
Problems,} VNU Science Press, Utrecht, 1987.
%
\bi{LZ82} B.\ M.\ Levitan and V.\ V.\ Zhikov, {\it Almost Periodic
Functions and Differential Equations}, Cambridge University Press,
Cambridge, 1982.
%
\bi{Ma92} A.\ I.\ Markushevich, {\it Introduction to the Classical
Theory of Abelian Functions}, Amer. Math. Soc, Providence, R.I., 1992.
%
\bi{McK79} H.\ P.\ McKean, {\it Integrable systems and algebraic
curves}, Global {A}nalysis, Lecture Notes in Math. {\bf
755}, Springer, Berlin, 1979, pp.\ 83--200.
%
\bi{MM75}  H.\ P.\ McKean and P.\ van Moerbeke, {\it The spectrum
of Hill's equation}, Invent. Math. {\bf 30}, 217--274 (1975).
%
%
\bi{MM80} H.\ P.\ McKean and P.\ van Moerbeke, {\it Hill and {T}oda
curves}, Comm. Pure Appl. Math. {\bf 33}, 23--42 (1980).
%
\bi{Mu77} D.\ Mumford, {\it An algebro-geometric construction of
commuting operators and of solutions to the {T}oda lattice
equation, {K}orteweg de{V}ries equation and related non-linear
equations}, Proceedings of the International Symposium on
Algebraic Geometry (Kyoto Univ., Kyoto, 1977), Kinokuniya
Book Store, Tokyo, 1978, pp.\ 115--153.
%
\bi{Mu84} D.\ Mumford, {\em Tata Lectures on Theta {II}}, Birkh\"
auser, Boston, 1984.
%
\bi{Na62} P.\ B.\ Na{\u\i}man, {\it On the theory of periodic and
limit-periodic Jacobian matrices}, Sov. Math. Dokl. {\bf 3}, 383--385
(1962).
%
\bi{Na64} P.\ B.\ Na{\u\i}man, {\it On the spectral theory of
non-symmetric periodic Jacobi matrices}, Zap. Meh.-Mat. Fak. Har{\c p}rime
kov. Gos. Univ. i Har{\c p}rime kov. Mat. Ob\v s\v c. (4) {\bf 30},
138--151 (1964). (Russian.)
%
\bibitem{NMPZ84} S.\ Novikov, S.\ V.\ Manakov, L.\ P.\ Pitaevskii, and
V.\ E.\ Zakharov, {\it Theory of Solitons}, Consultants Bureau,
New York, 1984.
%
\bi{Ro63} F.\ S.\ Rofe-Beketov, {\it The spectrum of
non-selfadjoint differential operators with periodic
coefficients}, Sov. Math. Dokl. {\bf 4}, 1563--1566 (1963).
%
\bibitem{RS78} M.\ Reed and B.\ Simon, {\it Methods of Modern Mathematical
Physics. IV: Analysis of Operators,} Academic Press, New York,
1978.
%
\bi{Sc65} G.\ Scharf, {\it Fastperiodische Potentiale}, Helv.
Phys. Acta, {\bf 38}, 573--605 (1965).
%
\bi{Se60} M.\ I.\ Serov, {\it Certain properties of the spectrum of a
non-selfadjoint differential operator of the second order}, Sov. Math.
Dokl. {\bf 1}, 190--192 (1960).
%
\bi{Si82} B.\ Simon, {\it Almost periodic Schr\"odinger operators:
A review}, Adv. Appl. Math. {\bf 3}, 463--490 (1982).
%
\bi{Te00} G.\ Teschl, {\it Jacobi Operators and Completely
Integrable Nonlinear Lattices}, Amer. Math. Soc., Providence, R.I., 2000.
%
\bi{To89} M.\ Toda, {\it Theory of Nonlinear Lattices},
Springer, Berlin, 1989.
%
\bi{To189} M.\ Toda, {\it {N}onlinear {W}aves and {S}olitons},
Kluwer, Dordrecht, 1989.
%
\bi{Mo79} P.\ van Moerbeke, {\it About isospectral deformations of
discrete Laplacians}, Global Analysis, Lecture Notes in Math. {\bf
755}, Berlin, 1979, pp.\ 313--370.
%
\bi{MM79} P.\ van Moerbeke and D.\ Mumford, {\it The spectrum of
difference operators and algebraic curves}, Acta Math. {\bf 143},
93--154 (1979).
%

\end{thebibliography}
\end{document}